\documentclass[aop,seceqn,MSNbibl,citesort,dvips]{arximspdf}

\doi{10.1214/09-AOP498}
\volume{38}
\issue{2}
\pubyear{2010}
\firstpage{896}
\lastpage{926}

\makeatletter

\newtheorem{theorem}{Theorem}[section]
\newtheorem{lemma}[theorem]{Lemma}
\newtheorem{proposition}[theorem]{Proposition}
\newtheorem{corollary}[theorem]{Corollary}

\makeatother

\begin{document}
\begin{frontmatter}

\title{Thick points of the Gaussian free field}
\runtitle{Thick points of the Gaussian free field}

\begin{aug}
\author[A]{\fnms{Xiaoyu} \snm{Hu}\thanksref{t1}\ead[label=e1]{xyhu@gucas.ac.cn}},
\author[B]{\fnms{Jason} \snm{Miller}\corref{}\thanksref{t2}\ead[label=e2]{jmiller@math.stanford.edu}} and
\author[C]{\fnms{Yuval} \snm{Peres}\ead[label=e3]{peres@microsoft.com}}
\runauthor{X. Hu, J. Miller and Y. Peres}
\affiliation{Chinese Academy of Sciences, Stanford University and
Microsoft Research}
\address[A]{X. Hu\\
Chinese Academy of Sciences\\
Department of Mathematics\\
Beijing 100049\\
China\\
\printead{e1}}
\address[B]{J. Miller\\
Stanford University\\
Department of Mathematics\\
Stanford, California 94305\\
USA\\
\printead{e2}}
\address[C]{Y. Peres\\
Microsoft Research\\
Theory Group\\
Redmond, Washington 98052\\
USA\\
\printead{e3}}
\end{aug}

\thankstext{t1}{Supported in part by NNSF of China Grant 10871200.}
\thankstext{t2}{Supported in part by NSF Grants DMS-04-06042 and DMS-08-06211.}

\received{\smonth{3} \syear{2009}}
\revised{\smonth{7} \syear{2009}}

%
\begin{abstract}
Let $U \subseteq\mathbf{C}$ be a bounded domain with smooth boundary
and let $F$ be an instance of the continuum Gaussian free field on $U$
with respect to the Dirichlet inner product $\int_U \nabla f(x) \cdot
\nabla g(x) \,dx$. The set $T(a;U)$ of $a$-thick points of $F$ consists
of those $z \in U$ such that the average of $F$ on a disk of radius $r$
centered at $z$ has growth $\sqrt{a/\pi} \log\frac{1}{r}$ as $r \to
0$. We show that for each $0 \leq a \leq2$ the Hausdorff dimension of
$T(a;U)$ is almost surely $2-a$, that $\nu_{2-a}(T(a;U)) = \infty$ when
$0 < a \leq2$ and $\nu_{2}(T(0;U)) = \nu_2(U)$ almost surely, where
$\nu_\alpha$ is the Hausdorff-$\alpha$ measure, and that $T(a;U)$ is
almost surely empty when $a > 2$. Furthermore, we prove that $T(a;U)$
is invariant under conformal transformations in an appropriate sense.
The notion of a thick point is connected to the Liouville quantum
gravity measure with parameter $\gamma$ given formally by $\Gamma(dz) =
e^{\sqrt{2\pi} \gamma F(z)} \,dz$ considered by Duplantier and Sheffield.
\end{abstract}

%

\end{frontmatter}

\section{Introduction}

Let $U \subseteq\mathbf{C}$ be a bounded domain with smooth boundary
and for
$f,g \in C_0^\infty(U)$ let
\[
(f,g)_\nabla= \int_U \nabla f (x) \cdot\nabla g(x)\,dx
\]
denote the Dirichlet inner product of $f$ and $g$ where $dx$ is the
Lebesgue measure. Let $(f_n)$ be an orthonormal basis of the Hilbert
space closure $H_0^1(U)$ of $C_0^\infty(U)$ under $(\cdot,\cdot
)_\nabla
$. The continuum Gaussian free field (GFF) $F=F_U$ on $U$ is given
formally as a random linear combination
%
\begin{equation}
\label{eqn::gff_definition}
F = \sum_n \alpha_n f_n,
\end{equation}
where $(\alpha_n)$ is an i.i.d. Gaussian sequence.

The GFF is a $2$-time dimensional analog of the Brownian motion. Just
as the Brownian motion can be realized as the scaling limit of many
random curve ensembles, the GFF arises as the scaling limit of a number
of random surface ensembles \cite{BAD96,KEN01,NS97,RV08}.
The purpose of this article is to study the fractal
geometry and conformal invariance properties of its extremal points. It
is not possible to make sense of $F$ as a function since the sum in
(\ref{eqn::gff_definition}) does not converge in a topology that would
allow us to do so. However, it does converge almost surely in the space
of distributions and is sufficiently regular that there is no
difficulty in interpreting its integral with respect to Lebesgue
measure over sufficiently nice Borel sets. This class includes, for
example, disks, squares and the conformal images of such. Thus, to make
the notion of an extremal point precise, we first regularize by
averaging the field over disks of radius $r$ and then study those
points where the average is unusually large as $r \to0$.

With this is in mind, we say that $z$ is an \textit{$a$-thick point} provided
%
\begin{equation}
\label{eqn::thick_def}
\lim_{r \to0} \frac{\mu(D(z,r))}{\pi r^2 \log 1/r} =
\sqrt{\frac{a}{\pi}},
\end{equation}
where $D(z,r)$ denotes the disk of radius $r$ centered at $z \in U$ and
$\mu(A) = \int_{A} F(x) \,dx$. While integration against the GFF does not
define a measure, we can still interpret the quantity $\mu(A)$ as
measuring the signed mass the GFF associates with $A$, whenever it is
defined. This motivates our usage of the term ``thick point,'' which
has become standard terminology in the literature which studies the
extremes of the occupation measure of a stochastic process
\cite{DAV05,DPRZ99,DPRZ00,DPRZ01}.

Let $T(a;U)$ denote the set of $a$-thick points of $F$ and let $\nu
_\alpha$ denote the Hausdorff-$\alpha$ measure.
\begin{theorem}
\label{thm::dimension}
For any $0 \leq a \leq2$, the Hausdorff dimension of $T(a;U)$ is
almost surely $2-a$. Moreover,
\[
\mathbf{P}[\nu_2(T(0;U)) = \nu_2(U)] = 1 \quad\mbox{and}\quad \mathbf
{P}[\nu_{2-a}(T(a;U)) =\infty] = 1,
\]
when $0 < a \leq2$. In particular, $\mathbf{P}[|T(2;U)| = \infty] = 1$.
Finally, if $a > 2$ then $T(a;U)$ is almost surely empty.
\end{theorem}

The proof of this theorem easily extends to the five other cases where
one replaces the limit in (\ref{eqn::thick_def}) with either $\liminf$
or $\limsup$ and the equality with ``not less than.'' The particular
choice of normalization in (\ref{eqn::thick_def}) is so that the
dimension is a linear function of $a$.

It is possible to make sense of the \textit{circle average process}
\[
F(z,r) = \frac{1}{2\pi r}\int_{\partial D(z,r)} F(x) \sigma(dx),
\]
$\sigma(dx)$ the length measure, in such a way that it is a continuous
function in $(z,r)$ \cite{DS08,SHE06}. We will describe this
construction in the next section and, furthermore, argue that almost
surely
%
\begin{equation}
\label{eqn::polar}
\int_0^r 2\pi s F(z,s) \,ds = \int_{D(z,r)} F(x) \,dx \qquad\mbox{for all } (z,r).
\end{equation}
This gives rise to another collection of thick points, namely the set
$T^C(a;U)$ consisting of those $z \in U$ satisfying
\[
\lim_{r \to0} \frac{1}{\log 1/r} F(z,r) = \sqrt{\frac
{a}{\pi}}.
\]
Our proof implies that the Hausdorff dimension of $T^C(a;U)$ is $2-a$
almost surely, and we include this result as a separate theorem.
\begin{theorem}
\label{thm::dimension_circle}
For any $0 \leq a \leq2$, the Hausdorff dimension of $T^C(a;U)$ is
almost surely $2-a$. Moreover,
\[
\mathbf{P}[\nu_2(T^C(0;U)) = \nu_2(U)] = 1 \quad\mbox{and}\quad \mathbf
{P}[\nu_{2-a}(T^C(a;U))
= \infty] = 1,
\]
when $0 < a \leq2$. In particular, $\mathbf{P}[|T^C(2;U)| = \infty]
= 1$.
Finally, if $a > 2$, then $T^C(a;U)$ is almost surely empty.
\end{theorem}

As before, our proof also extends to the cases where one replaces the
limit with either a $\liminf$ or $\limsup$ and the equality with ``not
less than.''

Suppose that $V$ is another domain, $\varphi\dvtx U \to V$ is a
conformal transformation, and for $A \subseteq V$, formally set
%
\begin{equation}
\label{eqn::gff_on_v}
F_V = F_U \circ\varphi^{-1} \quad\mbox{and}\quad \mu_V(A) = \int_A F_V(x) \,dx.
\end{equation}
As the Dirichlet inner product is invariant under precomposition by
conformal maps, it follows that $F_V$ has the law of a GFF on $V$. Our
next theorem is a uniform estimate on the difference between $\mu
_U(D(\xi,r))$ and $\mu_V(D(\varphi(\xi),|\varphi'(\xi
)|r))\times|\varphi'(\xi)|^{-2}$.
\begin{theorem}
\label{thm::conf_invariance}
If $K \subseteq U$ is compact, then almost surely
%
\begin{eqnarray}
\label{eqn::conf_error_estimate}
&&\lim_{r \to0} \sup_{\xi\in K} \frac{1}{\pi r^2 \log 1/r}
\bigl|\mu_{U}(D(\xi,r)) \nonumber\\[-8pt]\\[-8pt]
&&\hspace*{91.5pt}{} - \mu_V(D(\varphi(\xi),|\varphi'(\xi
)|r))|\varphi'(\xi)|^{-2} \bigr| = 0.\nonumber
\end{eqnarray}
\end{theorem}

An immediate consequence of this is the conformal invariance of the
thick points.
\begin{corollary}
\label{cor::conf_invariance}
The set of thick points is a conformal invariant. More precisely, if
$T(a;V)$ denotes the $a$-thick points of $\mu_V$ as in (\ref
{eqn::gff_on_v}), then
\[
\mathbf{P}\bigl(\varphi(T(a;U)) = T(a;V) \mbox{ for all } 0 \leq a \leq
2\bigr) = 1.
\]
\end{corollary}

Let $G = (V,E)$ be a finite graph with distinguished subset $V_\partial
\subseteq V$. The law of the discrete GFF (DGFF) is given by the Gibbs
measure with Hamiltonian $H(f) = \frac{1}{2} \sum_{x \sim y} (f(x) -
f(y))^2$ for $f|V_\partial\equiv0$. The Ginzburg--Landau (GL) $\nabla
\phi$ interface model, also known as the anharmonic crystal, is a
non-Gaussian analog of the DGFF and arises by replacing
\mbox{$|\cdot|^2$} in $H(f)$ with a symmetric, convex function with
quadratic growth. The behavior of the extremal points of the DGFF and
GL model are studied in \cite{BDG01,DAV06,DG00} in the special case
that $V$ is a lattice approximation of a smooth subset in $\mathbf{C}$.
In Theorem 1.3 of \cite{DAV06}, Daviaud shows that if for each
$\varepsilon> 0$ one lets $F_\varepsilon$ have the law of the DGFF on
the induced subgraph $U_\varepsilon= U \cap\varepsilon\mathbf{Z}^2$
then the cardinality of the set $\mathcal{H} _\varepsilon(a) = \{z \in
U_\varepsilon\dvtx F_\varepsilon(x) \geq\sqrt {a/\pi}
\log\frac{1}{\varepsilon}\}$ of ``\mbox{$a$-high} points'' has growth
$\varepsilon ^{-(2-a)}$ as $\varepsilon\to0$. This growth exponent
represents a sort of discrete Hausdorff dimension so that this is the
natural discrete analog of Theorem~\ref{thm::dimension}. An interesting
open question is to see if this result for the DGFF or analogous
results for the models considered in \cite{KEN01,NS97,RV08} can be
deduced from Theorem \ref{thm::dimension}. A natural starting point to
answering this question is the coupling of the DGFF and the GFF
suggested in \cite{SHE06}. A proof of the reverse implication seems
less hopeful since, intuitively, a point $z$ is $a$-thick if and only
if it is an $a$-high point for $F_\varepsilon$, all sufficiently small
$\varepsilon> 0$ and the result of \cite{DAV06} only provides estimates
of the number of $a$-high points for a fixed $\varepsilon> 0$.

Suppose that $U$ is simply connected. If $(S,g)$ is a Riemann surface
homeomorphic to $U$, then the classical uniformization theorem implies
that $(S,g)$ is conformally equivalent either to $U$ or $\mathbf{C}$. The
former case is in turn equivalent to the existence of global
coordinates with respect to which the metric $g$ takes the form
$e^{\lambda(z)} \,dz$ for some $\lambda\dvtx U \to\mathbf{R}$. One natural
construction of a random surface is to select $\lambda\dvtx U \to
\mathbf{R}$
randomly and then take the surface with metric $e^{\lambda(z)}\,dz$. Fix
$0 \leq\gamma< 2$. Formally, the Liouville quantum gravity with
parameter $\gamma$ corresponds to the case $\lambda(z) = \sqrt{2\pi}
\gamma F(z)$. This, however, does not make sense since $F$ is not even
pointwise defined. To make this rigorous, Duplantier and Sheffield in
\cite{DS08} consider the random surfaces with continuous metric
$r^{\gamma^2/2} e^{\sqrt{2\pi} \gamma F(z,r)} \,dz$ and study their
behavior as $r \to0$. Although understanding the limiting object as a
metric space is still out of reach, they show that the associated
random area measures $\Gamma_r$ have a weak limit $\Gamma$ as $r \to
0$. For $A$ Borel, the quantity $\Gamma(A)$ is referred to as the
$\gamma$-quantum area of $A$. It is shown in Proposition 3.4 of \cite
{DS08} and the discussion thereafter that $\Gamma$ is almost surely
supported on
\[
T_{\geq}^{C,i}(a;U) = \Biggl\{ z \in U \dvtx\liminf_{r \to0} \frac
{1}{\log
1/r} F(z,r) \geq\sqrt{\frac{a}{\pi}} \Biggr\},
\]
where $a = \gamma^2/2$.
Thus, developing an understanding of the geometry of the thick points
leads to an understanding of the geometry of the support of $\Gamma$.
Note that our definition is slightly different from that appearing in
\cite{DS08} because of a difference in a choice of normalization of the
Dirichlet inner product. Denote by $\tilde{D}(z,r) = \sup\{s \dvtx
\Gamma
(D(z,s)) \leq r\}$ the quantum ball of radius $r$ centered at $z$. Let
$X \subseteq U$ be a random Borel set independent of $F$ and let $X^r =
\bigcup_{z \in X} \tilde{D}(z,r)$ be the $r$-quantum neighborhood of $X$.
Then $X$ is said to have quantum scaling expectation exponent $\Delta$ provided
\[
\lim_{r \to0} \frac{\mathbf{E}\Gamma(X^r)}{\log r} = \Delta.
\]
Duplantier and Sheffield speculate (\cite{DS08}, page 26) that if $X$ has
quantum expectation scaling exponent $\Delta$, then its quantum support
is concentrated on $T_\geq^{C,i}(\alpha;U)$ where $\alpha= (\gamma-
\gamma\Delta)^2/2$.

The remainder of the paper is organized as follows. In Section \ref
{sec::gff_background}, we will give a brief overview of the basic
properties of the GFF; see \cite{SHE06} for a more thorough
introduction and \cite{JAN97} for more on the closely related notion of
a Gaussian Hilbert space. Next, in Section \ref{sec::hd}, we prove
Theorems \ref{thm::dimension} and \ref{thm::dimension_circle}. The
first step is to establish the identity (\ref{eqn::polar}) which is a
consequence of the fact that the Lebesgue measure on a disk can be
written as a limit of Riemann sums of the length measure where the
convergence is an appropriate Sobolev space. This allows us to sandwich
the sets considered in Theorems \ref{thm::dimension} and \ref
{thm::dimension_circle} between $T^C(a;U)$ and
\[
T_\geq^{C,s}(a;U) = \Biggl\{ \limsup_{r \to0} \frac{1}{\log
1/r} F(z,r) \geq\sqrt{\frac{a}{\pi}} \Biggr\}
\]
so that we need only show $\dim_H T^C(a;U) \geq2-a$ and $\dim_H
T_\geq
^{C,s}(a;U) \leq2-a$. We prove the more difficult lower bound using a
multi-scale refinement of the second moment method, similar to the
techniques employed in \cite{DPRZ01}. Roughly speaking, the crucial
estimate that one needs is a quantitative bound on the degree to which
the events that two given points are $a$-thick are approximately
independent. We address this by considering a special subset which we
term ``perfect thick points.'' These are defined in such a way so that
the approximate independence is a consequence of the Markov property of
the field. The upper bound follows from an estimate of the modulus of
continuity of $F(z,r)$ and that for fixed $z$ the processes $r \mapsto
F(z,e^{-r})$ evolve as Brownian motions.

Finally, in Section \ref{sec::conf_invariance}, we prove Theorem \ref
{thm::conf_invariance}. It is easy to predict that the conformal
invariance result is true since the GFF itself is conformally
invariant, thick points are defined in terms of averages over small
disks, and conformal maps send disks to disks at infinitesimal scales.
This intuition, however, is far from a proof since integration against
the GFF does not define a measure, much less a measure that is
absolutely continuous with respect to Lebesgue measure. In particular,
the GFF can assign large mass to a set that is small in the Lebesgue
sense precisely due to the presence of thick points. The basic idea of
our proof is as follows. We use the Markov property of the field to
reduce to the case that $V = [0,1]^2$. This choice is particularly
convenient because the $H_0^1([0,1]^2)$ orthonormal basis $(f_n)$ given
by the eigenvectors of the Laplacian is given by products of sine
functions; this makes many of our computations elementary and explicit.
We then show that if $A$ is a small dyadic square centered at $z$ or
the image of such centered at $\xi= \varphi^{-1}(z)$ under $\varphi$
then then $\mu_V(A)$ is sufficiently well approximated by $\sum_{n=1}^N
\alpha_n f_n(z) |A|$. Using a covering argument, we then deduce that an
analogous estimate also holds for disks $D(z,r)$ and conformal images
of disks $\varphi(D(\xi,r))$ with small radii. This argument is
sensitive to the geometry of a disk since we need that the number
$N(t;D(z,r))$ of maximal dyadic squares of side length $t$ in $D(z,r)$
does not grow too quickly as $t \to0$. Theorem \ref
{thm::conf_invariance} then follows from a bound on the Lebesgue
measure of the symmetric difference $\varphi(D(\xi,r)) \Delta
D(\varphi
(\xi),|\varphi'(\xi)|r)$ and some Gaussian estimates.

Throughout the paper, we will make use of the following notation. If
$f,g$ are two functions, then we write $f \sim g$ provided that
$f(t)/g(t) \to1$ as either $t \to\infty$ or $t \to0$, the case being
clear from the context. If $f_\alpha, g_\alpha$ are one-parameter
families of functions, then $f_\alpha\sim g_\alpha$ uniformly means
that $f_\alpha(t)/g_\alpha(t) \to1$ uniformly in $\alpha$. We say that
$f = O(g)$ if there exists a constant $C > 0$ such that $|f(t)| \leq
C|g(t)|$ for all $t$ and that $f = o(g)$ provided that $|f(t)|/|g(t)|
\to0$ as either $t \to0$ or $t \to\infty$, the case being clear from
the context. Finally, we say $f_\alpha= O(g_\alpha)$ and $f_\alpha=
o(g_\alpha)$ uniformly in $\alpha$ if the constant and convergence are
uniform in $\alpha$, respectively.

\section{The Gaussian free field}
\label{sec::gff_background}

The purpose of this section is to recall the basic properties of the
GFF. Let $U$ be a bounded domain in $\mathbf{C}$ with smooth boundary
and let
$C_0^\infty(U)$ denote the set of $C^\infty$ functions compactly
supported in $U$. We begin with a short discussion of Sobolev spaces;
the reader is referred to Chapter 5 of \cite{EVAN02} or Chapter 4 of
\cite{TAY96} for a more thorough introduction. With $\mathbf{N}_0 = \{
0,1,\ldots
\}$, the nonnegative integers when $f \in C_0^\infty(U)$ and $\alpha=
(\alpha_1,\alpha_2) \in\mathbf{N}_0^2$ we let $D^\alpha f =
\partial
_1^{\alpha
_1} \partial_2^{\alpha_2} f$. For $k \in\mathbf{N}_0$, we define
the $H^k(U)$-norm
%
\begin{equation}
\label{eqn::sobolev}
\|f\|_{H^k(U)}^2 = \sum_{|\alpha| \leq k} \int_U |D^\alpha f(x)|^2
\,dx,
\end{equation}
where $|\alpha| = \alpha_1 + \alpha_2$. The Sobolev space $H_0^k(U)$ is
given by the Banach space closure of $C_0^\infty(U)$ under \mbox{$\Vert
\cdot\Vert_{H^k(U)}$}. If $s \geq0$ is not necessarily an integer, then
$H_0^s(U)$ can be constructed via the complex interpolation of
$H_0^0(U) = L^2(U)$ and $H_0^k(U)$ where $k \geq s$ is any positive
integer (see Chapter 4 Section 2 of \cite{TAY96} for more on this
construction and also Chapter 4 of \cite{KAT04} for more on
interpolation). A consequence of this is that if $T \dvtx C_0^\infty
(U) \to C_0^\infty(U)$ is a linear map continuous with respect to the
$L^2(U)$ and $H^k(U)$ topologies, then it is also continuous with
respect to $H^s(U)$ for all $0 \leq s \leq k$. For $s \geq0$, we
define $H^{-s}(U)$ to be the Banach space dual of $H_0^s(U)$ where the
dual pairing of $f \in H^{-s}(U)$ and $g \in H_0^s(U)$ is given
formally by the usual $L^2(U)$ inner product
\[
(f,g) = (f,g)_{L^2(U)} = \int_U f(x) g(x) \,dx.
\]
More generally, for any $s \in\mathbf{R}$ the $H^s(U)$-topology can be
constructed explicitly via the norm
%
\begin{equation}
\label{eqn::fourier_sobolev}
\|f\|_s^2 = \int(1+\xi_1^2 + \xi_2^2)^s (\widehat{f}(\xi))^2\,d\xi,
\end{equation}
where
\[
\widehat{f}(\xi) = \int e^{-i \xi\cdot x} f(x) \,dx
\]
is the Fourier transform of $f$. We will be most interested in the
space $H_0^1(U)$. An application of the Poincare inequality (Chapter 4,
Proposition 5.2) gives that the norm induced by the Dirichlet inner product
\[
(f,g)_\nabla= \int_U \nabla f \cdot\nabla g \qquad\mbox{for } f,g \in
C_0^\infty(U),
\]
is equivalent to \mbox{$\| \cdot\|_{H^1(U)}$}. This choice of inner product
is particularly convenient because it is invariant under precomposition
by conformal transformations.

The GFF $F=F_U$ on $U$ is given formally as a random linear combination
of an orthonormal basis $(f_n)$ of $H_0^1(U)$
\[
F = \sum_n \alpha_n f_n,
\]
where $(\alpha_n)$ is an i.i.d. sequence of standard Gaussian. Although
the sum defining $F$ does not converge in $H_0^1(U)$, for each
$\varepsilon> 0$ it does converge almost surely in $H^{-\varepsilon}(U)$
(\cite{SHE06}, Proposition 2.7 and the discussion thereafter) and, in
particular, $H^{-1}(U)$. If $f,g \in C_0^\infty(U)$, then an
integration by parts gives $(f,g)_\nabla= -( f,\Delta g)$. Using this,
we define
\[
(F,f)_\nabla= -(F, \Delta f) \qquad\mbox{for } f \in C_0^\infty(U).
\]
Observe that $(F,f)_\nabla$ is a Gaussian random variable with mean
zero and variance $(f,f)_\nabla$. Hence, by polarization, $F$ induces a
map $C_0^\infty(U) \to\mathcal{G}$, $\mathcal{G}$ a Gaussian
Hilbert space, that
preserves the Dirichlet inner product. This map extends uniquely to
$H_0^1(U)$ which allows us to make sense of $(F,f)_\nabla$ for all $f
\in H_0^1(U)$. We are careful to point out, however, that while $(F,
\cdot)_\nabla$ is well defined off of a set of measure zero as a linear
functional on $C_0^\infty(U)$ this is not the case for general $f \in
H_0^1(U)$. This is a technical point that we will touch on a bit later.
It is not hard to see that the law of $F$ is independent of the choice
of $(f_n)$; the eigenvectors of the Laplacian serve as a convenient
choice since they are also orthogonal in $L^2(U)$. In particular, when
$U = [0,1]^2$ then $F_{[0,1]^2}$ admits the explicit representation
%
\begin{equation}
\label{eqn::gff_square_rep}
F_{[0,1]^2}(x,y) = \sum_{i,j \geq1} \frac{2 \alpha_{ij}}{\pi\sqrt
{i^2 + j^2}} \sin(\pi i x) \sin(\pi j y) \qquad\mbox{for } (x,y) \in[0,1]^2.\hspace*{-30pt}
\end{equation}
If $V \subseteq\mathbf{C}$ is another domain, $\varphi\dvtx U \to
V$ is a
conformal transformation and $f,g \in C_0^\infty(U)$, then a change of
variables shows that the Dirichlet inner product is invariant under
precomposition by $\varphi^{-1}$:
\[
\int_V \nabla(f \circ\varphi^{-1}) \cdot\nabla(g \circ\varphi^{-1})
= (f,g)_\nabla.
\]
Thus if $(f_n)$ is an orthonormal basis of $H_0^1(U)$, then $(f_n \circ
\varphi^{-1})$ is an orthonormal basis of $H_0^1(V)$, so that if $F$ is
a GFF on $U$, then $F_V = F \circ\varphi^{-1}$ has the law of a GFF
on $V$.

If $\eta\in H^{-1}(U)$ so that $-\Delta^{-1} \eta\in H_0^1(U)$, then
$(F,-\Delta^{-1} \eta)_\nabla= ( F, \eta)$. The particular case that
will be of interest to us is when $\eta(z,r)$ is the uniform measure on
the circle $\partial D(z,r)$ where we think of $F(z,r) = (F,\eta(z,r))$
as the mean value of $F$ on $\partial D(z,r)$. Letting
\[
G(x,y) = -\frac{1}{2\pi} \bigl({\log}|x-y| -\phi^y(x) \bigr),
\]
where for each fixed $y \in U$ we denote by $x \mapsto\phi^y(x)$ the
harmonic extension of ${\log}|x-y|$ from $\partial U$ to $U$, be the
Green's function for the Dirichlet problem of the Laplacian on $U$ with
zero boundary conditions, observe
\begin{eqnarray*}
\mathbf{E}( F, - \Delta^{-1} f )_\nabla( F, -\Delta^{-1} g )_\nabla
&=& -( f, \Delta^{-1} g )\\
&=& \int_U \int_U f(x) g(y) G(x,y) \,dx\, dy.
\end{eqnarray*}
When $D(z,e^{-t_1}) \subseteq U$ and $s,t > t_1$ we have
\begin{eqnarray*}
\mathbf{E}F(z,e^{-s}) F(z,e^{-t})
&=& \mathbf{E}(F, - \Delta^{-1} \eta(z,e^{-s}) )_\nabla( F, -\Delta^{-1}
\eta
(z,e^{-t}) )_\nabla\\
&=& \frac{s}{2\pi} + C(z),
\end{eqnarray*}
where $C(z)$ is a constant depending only on $z$ and not $s,t$. Hence,
with $z \in U$ fixed and letting $B(z,t) = \sqrt{2\pi} F(z,e^{-t})$ the
process $t \mapsto B(z,t) - B(z,t_1)$ is Gaussian with the mean and
autocovariance of a standard Brownian motion.

Using the Kolmogorov--Centsov theorem one can show (\cite{DS08},
Proposition 3.1) that $(z,r) \mapsto F(z,r)$ has a locally $\gamma$-H\"
{o}lder continuous modification whenever $\gamma< 1/2$ is fixed. We
will need some control of the H\"{o}lder norm of $F(z,r)$ on compact
intervals as $r \to0$; we are able to do this using Lemma \ref
{lem::kolm_centsov}, a refinement of the Kolmogorov--Centsov theorem.
\begin{proposition}
\label{prop::holder}
The circle average process $F(z,r)$ possesses a modification $\tilde
{F}(z,r)$ such that for every $0 < \gamma< 1/2$ and $\varepsilon,\zeta>
0$ there exists $M = M(\gamma,\varepsilon,\zeta)$ such that
%
\begin{equation}
\label{eqn::holder}
|\tilde{F}(z,r) - \tilde{F}(w,s)| \leq M \biggl( \log\frac
{1}{r}
\biggr)^\zeta\frac{|(z,r)-(w,s)|^\gamma}{r^{\gamma+\varepsilon}}
\end{equation}
for all $z,w \in U$ and $r,s \in(0,1]$ with $1/2 \leq r/s \leq2$.
\end{proposition}
\begin{pf}
Note that if $\tilde{F}$ and $\tilde{F}'$ are two different
modifications satisfying (\ref{eqn::holder}), then they are almost
surely equal by continuity. Thus, it suffices to show that $F$
satisfies the hypotheses of Lemma \ref{lem::kolm_centsov} for $\alpha
,\beta$ arbitrarily large with $\beta/\alpha$ arbitrarily close to
$1/2$. With $a = (z,w,r,s) \in\mathcal{U}= U^2 \times[0,\infty)^2$,
we know that
\begin{eqnarray*}
\Upsilon(a)
&=& \mathbf{E}F(z,r) F(w,s) = (\eta(z,r), -\Delta^{-1} \eta(w,s))\\
&=& ( -\Delta^{-1} \eta(z,r), \eta(w,s)).
\end{eqnarray*}
One can check directly (see the discussion after Proposition 3.1 of
\cite{DS08}) that $\xi_r^z = -\Delta^{-1} \eta(z,r)$ is given by
\[
\xi_r^z(y) = \tau_r^z(y) - \psi_r^z (y),
\]
where $\tau_r^z(y) = -\log\max(r,|z-y|)$ and $\psi_{r}^z$ is the
harmonic extension $\tau_r^z$ from $\partial U$ to $U$. As $|{\log}
\frac
{x}{y}| \leq\frac{|x-y|}{x \wedge y}$ for $x,y > 0$, we have
%
\begin{equation}
\label{eqn::tau_bound}
|\tau_r^z(y) - \tau_{r'}^{z'}(y')|
\leq C \frac{|r-r'| + |z-z'| + |y-y'|}{r \wedge r'}.
\end{equation}
In particular, this holds when $y_0 \in\partial U$. This implies that
the partial derivatives $\partial_y,\partial_z,\partial_r$ of $\psi
_r^z(y_0)$ are all $O(1/r)$ uniformly $z \in U$ and $y_0 \in\partial
U$ when $r \in(0,1]$. Since these partials are harmonic from the
maximum principle we conclude that (\ref{eqn::tau_bound}) holds with
$\psi_r^z$ in place of $\tau_r^z$. This gives
\[
|\xi_s^w(z + rx) - \xi_{s'}^{w'}(z' + r'x)|
\leq\frac{C |a-a'|}{s \wedge s'}
\]
for all $a,a' \in\mathcal{U}$ and $x \in\mathbf{S}^1$ so that
\[
|\Upsilon(a) - \Upsilon(a')|
\leq\int_{\mathbf{S}^1} |\xi_s^w(z + r x) - \xi_{s'}^{w'}(z' + t'x)|
\sigma(dx)
\leq\frac{C|a-a'|}{s \wedge s'} .
\]
As everything is symmetric,
\[
|\Upsilon(a) - \Upsilon(a')| \leq\frac{C|a-a'|}{(r \wedge r') \vee(s
\wedge s')}.
\]
Hence,
\begin{eqnarray*}
&&\mathbf{E}\bigl( F(z,r) - F(w,s)\bigr)^2\\
&&\qquad\leq |\mathbf{E}(F(z,r))^2 - \mathbf{E}F(z,r) F(w,s)| +
|\mathbf{E}(F(w,s))^2 - \mathbf{E}F(z,r) F(w,s)|\\
&&\qquad= |\Upsilon(z,z,r,r) - \Upsilon(z,w,r,s)| + |\Upsilon(w,w,s,s) -
\Upsilon(z,w,r,s)|\\
&&\qquad\leq \frac{C|(z,r) - (w,s)|}{r \wedge s}.
\end{eqnarray*}
This implies that for any $\alpha> 1$, $z,w \in U$, and $r,s \in
(0,1]$ we have
\[
\mathbf{E}|F(z,r) - F(w,s)|^\alpha\leq C \biggl( \frac{|(z,r) - (w,s)|}{r
\wedge s} \biggr)^{\alpha/2},
\]
which puts us exactly in the setting of Lemma \ref{lem::kolm_centsov}.
\end{pf}

From now on, we assume that $F(z,r)$ is a modification as in
Proposition \ref{prop::holder}.

The most natural way to make sense of $\int_{A} F(x) \,dx$ is to show
that $1_A \in H_0^\varepsilon(U)$ for some $\varepsilon> 0$ and then to
interpret the integral as the dual pairing of $F$ and $f = 1_A$. To
show\vspace*{1pt} $\| 1_A \|_{H^\varepsilon(U)} < \infty$ it suffices to show that the
asymptotics of the Fourier transform\vspace*{1pt} $\widehat{f}$ are sufficiently
well-behaved so that the \mbox{$\| \cdot\|_\varepsilon$} norm of (\ref
{eqn::fourier_sobolev}) is finite. When $A$ is a disk or square then the
Fourier transform $\widehat{f}(r,\theta)$ of $1_A$ in polar coordinates
satisfies
\[
\sup_{r \geq0} r^{3/2-\varepsilon} |\widehat{f}(r,\theta)| \in
L^2(\mathbf{S}^1)
\]
for every $\varepsilon> 0$ (Theorems 1 and 2 of \cite{RAN69}). This
implies $1_A \in H_0^\varepsilon(U)$ whenever $\varepsilon< 1/2$. If
$\varphi\dvtx U \to V$ is a conformal transformation and $W \subseteq
U$ is an open set such that $\overline{W} \subseteq U$ then $c \leq
|\varphi'(z)| \leq C$ for all $z \in W$ and $0 < c \leq C < \infty$. It
thus follows that precomposition by $\varphi^{-1}$ induces a continuous
linear map $L^2(W) \to L^2(\varphi(W))$ and $H_0^1(W) \to
H_0^1(\varphi
(W))$. Therefore by interpolation $1_{\varphi(A)} = 1_A \circ\varphi
^{-1} \in H_0^\varepsilon(\varphi(W)) \subseteq H_0^\varepsilon(V)$ for all
$\varepsilon< 1/2$ when $A \subseteq W$ is a square or disk.

The Lebesgue measure $\rho= \rho(z,r)$ on $D(z,r)$ can be expressed as
the integral $\rho= \int_0^r 2\pi s \eta(z,s) \,ds$. This gives rise to
two different interpretations of $\int_{D(z,r)} F(x) \,dx$ both of which
will be important for us. The first is as the dual pairing we have
already mentioned and the second is
%
\begin{equation}
\label{eqn::disk_integral}
\int_0^r 2\pi s F(z,s) \,ds = \sqrt{2\pi} \int_{-\log r}^\infty B(z,t)
e^{-2t} \,dt.
\end{equation}
Thus, we must be careful to ensure that they agree in an appropriate
sense. This does not represent a serious difficulty, however, since it
is easy to see that the Riemann sums corresponding to $ \int_0^r 2\pi s
\eta(z,s) \,ds$ converge to $\rho(z,r)$ in $H^{-1}(U)$.
If $\Pi$ is any partition of $[0,r]$ then as random variables in
$\mathcal{G}$
\begin{eqnarray*}
\sqrt{2\pi}\biggl(F, \sum_{\Pi} t_k \eta(z,t_k)(t_{k+1} - t_k)\biggr)
&\stackrel{\mathrm{a.s.}}{=}& \sqrt{2\pi}\sum_{\Pi} (F, t_k \eta
(z,t_k))(t_{k+1} - t_k)\\
&\stackrel{\mathrm{a.s.}}{=}& \sum_{\Pi} B(z,{-\log t_k}) t_k (t_{k+1}
- t_k).
\end{eqnarray*}
Therefore,
\[
(F,\rho(z,r)) \stackrel{\mathrm{a.s.}}{=} \int_{{-}\log r}^\infty B(z,t)
e^{-2t} \,dt
\]
as random variables in $\mathcal{G}$. As both sides of the equation are
continuous in $(z,r)$, we obtain the following proposition.
\begin{proposition}
Almost surely,
\[
(F,\rho(z,r)) = \int_{{-}\log r}^\infty B(z,t) e^{-2t} \,dt\qquad \mbox{for all
} (z,r).
\]
In particular, $z$ is an $a$-thick point if and only if
%
\begin{equation}\qquad
\label{eqn::thick_point_bm}
\lim_{r \to0} \frac{\sqrt{2\pi} \int_{{-}\log r}^\infty B(z,t) e^{-2t}\,
dt}{\sqrt{\pi} r^2 \log 1/r}
= \lim_{r \to0} \frac{ \sqrt{2} \int_{{-}\log r}^\infty B(z,t) e^{-2t}\,
dt}{ r^2 \log 1/r}
= \sqrt{a}.
\end{equation}
\end{proposition}

Suppose that $W \subseteq U$ is an open set. Then there is a natural
inclusion of $H_0^1(W)$ into $H_0^1(U)$ given by the extension by value
zero. If $f \in C_0^\infty(W)$ and $g \in C_0^\infty(U)$, then as
$(f,g)_\nabla= -(f,\Delta g)$ it is easy to see that $H_0^1(U)$ admits
the orthogonal decomposition $\mathcal{M}\oplus\mathcal{N}$ where
$\mathcal{M}= H_0^1(W)$
and $\mathcal{N}$ is the set of functions in $H_0^1(U)$ that are
harmonic on
$W$. Thus, we can write
\[
F = F_W + H_W = \sum_n \alpha_n f_n + \sum_n \beta_n g_n,
\]
where $(\alpha_n),(\beta_n)$ are independent i.i.d. sequences of
standard Gaussians and $(f_n)$, $(g_n)$ are orthonormal bases of
$\mathcal{M}$
and $\mathcal{N}$, respectively. Observe that $F_W$ has the law of the
GFF on
$W$, $H_W$ the harmonic extension of $F|\partial W$ to $W$, and $F_W$
and $H_W$ are independent. We arrive at the following proposition.
\begin{proposition}[(Markov property)]
\label{prop::markov}
The conditional law of $F|W$ given $F | U \setminus W$ is that of the
GFF on $W$ plus the harmonic extension of the restriction of $F$ on
$\partial W$ to $W$. In particular, if $D(z,e^{-t_1}) \setminus
D(z,e^{-t_2})$ and $D(w,e^{-s_1}) \setminus D(w,e^{-s_2})$ are disjoint
annuli contained in $U$ then the Brownian motions $B(z,t) - B(z,t_1)$
for $t_1 \leq t \leq t_2$ and $B(w,s) - B(w,s_1)$ for $s_1 \leq s \leq
s_2$ are independent.
\end{proposition}

\section{The Hausdorff dimension}
\label{sec::hd}

Let $U$ be a bounded domain with smooth boundary. It follows from the
discussion in the previous section that we can express $T^C(a;U)$ and
$T_{\geq}^{C,s}(a;U)$ as
\begin{eqnarray*}
T^C(a;U) &=& \biggl\{z \in U \dvtx\lim_{t \to\infty} \frac{B(z,t)}{
\sqrt{
2}t} = \sqrt{a} \biggr\},\\
T_{\geq}^{C,s}(a;U) &=& \biggl\{z \in U \dvtx\limsup_{t \to\infty} \frac
{B(z,t)}{ \sqrt{ 2}t} \geq\sqrt{a} \biggr\}.
\end{eqnarray*}

\subsection{The upper bound}

\begin{lemma}
\label{lem::upper_bound}
If $0 \leq a \leq2$, then almost surely $\dim_H(T_{\geq}^{C,s}(a;U))
\leq2-a$. If $a > 2$, then $T_{\geq}^{C,s}(a;U)$ is empty.
\end{lemma}
\begin{pf}
First, we suppose that $0 \leq a \leq2$. Let $\varepsilon> 0$ be
arbitrary and take $K = \varepsilon^{-1}$. For each $n$, let $r_n =
n^{-K}$. With $\zeta\in(0,1)$, $\gamma\in(0,1/2)$ and $\tilde
{\gamma
} = (1+\varepsilon)\gamma$ fixed and $M = M(\gamma,\varepsilon,\zeta)$
as in
(\ref{eqn::holder}), we have
\begin{eqnarray*}
\biggl|B(z,t) - B\biggl(z,\log\frac{1}{r_n}\biggr) \biggr|
&=& \sqrt{2\pi} |F(z,e^{-t}) - F(z,r_n)|\\
&\leq& M K^\zeta(\log n)^\zeta\frac{(r_{n+1} - r_n)^\gamma
}{r_{n+1}^{\tilde{\gamma}}}\\
&=& O \bigl( (\log n)^\zeta n^{K \tilde{\gamma} - (K+1) \gamma} \bigr)\\
&=& O((\log n)^\zeta)
\end{eqnarray*}
uniformly in $n \in\mathbf{N}$, $z \in U$, and $\log\frac
{1}{r_{n}} < t
\leq
\log\frac{1}{r_{n+1}}$. Therefore, $z \in T_{\geq}^{C,s}(a;U)$ if and
only if
\[
\limsup_{n \to\infty} \frac{B(z, \log 1/r_n)}{\sqrt{2}
\log
1/{r_n}} \geq\sqrt{a}.
\]

For each $n \in\mathbf{N}$, let $(z_{nj})$ be a maximal
$r_n^{1+\varepsilon}$ net
of $U$. If $z \in D(z_{nj}, r_n)$, then
\[
\biggl|B\biggl(z, \log\frac{1}{r_n}\biggr) - B\biggl(z_{nj}, \log\frac
{1}{r_n}\biggr) \biggr|
\leq O((\log n)^\zeta).
\]
Let
\[
\delta(n) = C(\log n)^{\zeta-1} \quad\mbox{and}\quad \mathcal{I}_n =
\biggl\{ j \dvtx
\biggl|B\biggl(z_{nj},\log\frac{1}{r_n}\biggr)\biggr| \geq\sqrt{2}\bigl(\sqrt{a}- \delta(n)\bigr)
\log\frac{1}{r_n} \biggr\}.
\]
Then we see that for each $N \geq1$
\[
I(a,N) = \bigcup_{n \geq N} \{ D(z_{nj},r_n) \dvtx j \in\mathcal{I}_n\}
\]
is such that $z \in T_{\geq}^{C,s}(a;U)$ implies that there exists
arbitrarily small balls in $I(a,N)$ containing $z$ provided $C$ is
large enough.

Since $B(z,t)$ evolves as a Brownian motion Lemma \ref
{lem::norm_asymp_decay} implies
\[
\mathbf{P}(j \in\mathcal{I}_n)
= \mathbf{P} \Biggl( \frac{|B(z_{nj},\log 1/{r_n})|}{\sqrt
{\log
1/{r_n}}} \geq\bigl(\sqrt{a} - \delta(n)\bigr)\sqrt{2 \log\frac
{1}{r_n}} \Biggr)
= O \bigl( r_n^{a-o(1)} \bigr).
\]
Hence,
%
\begin{equation}
\label{eqn::cover_bound}
\mathbf{E}|\mathcal{I}_n| \leq O \biggl(\frac{r_n^{a-o(1)}}{ r_n^{2
(1+\varepsilon
)}}
\biggr) = O \bigl( r_n^{a-o(1) - 2 (1+\varepsilon)} \bigr).
\end{equation}
Letting $\alpha= 2-a+\frac{2 + a}{1+\varepsilon} \varepsilon$, we thus have
\begin{eqnarray*}
\mathbf{E} \biggl[\sum_{n \geq N} \sum_{j \in\mathcal{I}_n}
(\operatorname{diam}(D(z_{nj},r_n^{1+\varepsilon})))^\alpha\biggr]
&=& O \biggl(\sum_{n \geq N} r_n^{2\varepsilon- o(1)} \biggr)\\
&=& O \biggl( \sum_{n \geq N} n^{-2+o(1)} \biggr).
\end{eqnarray*}
This proves that the Hausdorff-$[2-a+ \frac{2 + a}{1+\varepsilon}
\varepsilon
]$ measure of $T_{\geq}^{C,s}(a;U)$ is $0$.

If $a > 2$, then all of our analysis still applies. In particular, for
$\varepsilon> 0$ such that $a > 2(1+\varepsilon)$ (\ref{eqn::cover_bound})
gives that $\mathbf{E}|\mathcal{I}_n| \to0$ as $n \to\infty$.
\end{pf}

\subsection{The lower bound}
\label{subsec::lower_bound}

Let $s_n = \frac{1}{ n!}$ and $t_n = -\log s_n$. Let $H \subseteq U$
be a fixed compact square. By rescaling, we may assume without loss of
generality that if $z \in H$, then $D(z,s_n) \subseteq D(z,s_1)
\subseteq U$ and that $H$ has side length $1$. We further assume $H =
[0,1]^2$ by translation. For $m \in\mathbf{N}$, let
\begin{eqnarray*}
E_{m}(z) &=& \Bigl\{ \sup_{t_m < t \leq t_{m+1}} \bigl|B(z,t) - B(z,t_m) - \sqrt
{2a} (t - t_m )\bigr| \leq\sqrt{t_{m+1}-t_m} \Bigr\},\\
F_m(z) &=& \Bigl\{ {\sup_{t \geq t_m}} |B(z,t) - B(z,t_m)| \leq(t-t_m) + 1\Bigr\}.
\end{eqnarray*}
We say that $z \in H$ is an $n$-perfect $a$-thick point provided that
the event $E^n(z) = \bigcap_{m \leq n} E_{m}(z) \cap F_{n+1}(z)$ occurs.
Note that on $E^n(z)$ for $t_m < t \leq t_{m+1}$ and $m \leq n$ we have
\begin{eqnarray*}
&&\bigl|B(z,t) - B(z,t_1) - \sqrt{2a}(t-t_1)\bigr|\\
&&\qquad\leq \sum_{k=1}^{m-1} \bigl|B(z,t_{k+1}) - B(z,t_{k}) - \sqrt
{2a}(t_{k+1}-t_k)\bigr| \\
&&\qquad\quad{} + \bigl|B(z,t) - B(z,t_{m}) - \sqrt{2a}(t-t_m)\bigr|\\
&&\qquad\leq\sum_{k=1}^{m} \sqrt{\log(k+1)} = o(m \log m) = o(t)
\qquad\mbox{as }m \to\infty,
\end{eqnarray*}
where we used $t_{n+1} - t_n = \log(n+1)$ in the last inequality.
Furthermore, if $t \geq t_{n+1}$ then
\[
|B(z,t) - B(z,t_1)| \leq\sum_{k=1}^{n} \sqrt{\log(k+1)} + O(t) = O(t).
\]
Divide $H$ into $s_n^{-2}$ squares of side length $s_n$. Let $C_n$
denote the set of centers of these squares and $C_n(a)$ the set of
centers in $H$ that are $n$-perfect. Finally, we let
\[
P(a) = \bigcap_{k \geq1} \overline{\bigcup_{n \geq k} \bigcup_{ z
\in
C_n(a)} S(z,s_n)}
\]
be the set of ``perfect $a$-thick points,'' with $S(z,r)$ denoting the
square centered at $z$ of side length $r$. We obtain the following
lemma as an immediate consequence of the continuity of $B(z,r)$.
\begin{lemma}
Almost surely $P(a) \subseteq T^C(a;U)$.
\end{lemma}
\begin{pf}
Fix $z \in P(a)$. Then there exists a sequence $(z_{n_k})_{k=1}^\infty$
so that $z_{n_k} \in C_{n_k}(a)$ for every $k$ and $|z_{n_k} - z| \leq
s_{n_k}$. Fix $t > 0$ and let $m$ be such that $t_m < t \leq t_{m+1}$.
Uniformly in $k$ such that $n_k > m+1$ we know that
\[
\bigl|B(z_{n_k},t) - B(z_{n_k},t_1) - \sqrt{2a}(t-t_1)\bigr| = o(t).
\]
Thus taking a limit as $k \to\infty$ and using the spatial continuity
of $B$, we have
\[
\bigl|B(z,t) - B(z,t_1) - \sqrt{2a}(t- t_1)\bigr| = o(t).
\]
Since $t > 0$ was arbitrary, dividing both sides by $t$ we arrive at
\[
\frac{|B(z,t) - \sqrt{2a} t|}{t} = o(1) \qquad\mbox{as } t \to\infty.
\]
Therefore, $z \in T^C(a;U)$.
\end{pf}

From now on, we let $\gamma_n = \prod_{k=1}^n \exp(\frac
{1}{2}\sqrt
{\log k})$.
\begin{lemma}
Suppose $z,w \in H$. Let $l \in\mathbf{N}$ be such that $w \in S(z,s_l)
\setminus S(z,s_{l+1})$. There exists $C > 0$ such that for every $n
\geq l$ we have
\[
\mathbf{P}\bigl(E^n(z) \cap E^n(w)\bigr) \leq C^l \gamma_l^{-a} s_l^{-a}
\mathbf{P}(E^n(z)) \mathbf{P}
(E^n(w)).
\]
\end{lemma}
\begin{pf}
By making the constant sufficiently large the inequality holds
uniformly when $l \leq2$, hence we assume that $l > 2$. Observe that
the events $E_{i}(z), E_{j}(w)$ for $l+1 < i \leq n$ and $1 \leq j \leq
n$, $j \neq l-1,l,l+1$ are independent. By adjusting $C > 0$ if
necessary, Lemma \ref{lem::bm_bound} gives us the bound
\[
\mathbf{P}(E_{m}(z)),\qquad \mathbf{P}(E_m(w)) \geq C \frac{\exp(
{a}/{2} \sqrt
{\log
m})}{m^a} \qquad\mbox{for all } 1 \leq m \leq n.
\]
By further adjusting $C > 0$, it follows from Proposition \ref
{prop::markov} that
\[
\mathbf{P} \biggl(\bigcap_{1 \leq i \leq l+1} E_{i}(z) \biggr)\mathbf
{P} \biggl(
\bigcap
_{l-1 \leq j \leq l+1}E_{j}(w) \biggr)
\geq C^l \gamma_l^{a} s_l^{a}.
\]
Applying Proposition \ref{prop::markov} again, we therefore have the inequality
\begin{eqnarray*}
&&\mathbf{P}\bigl(E^n(z) \cap E^n(w)\bigr)
\\
&&\qquad\leq \mathbf{P} \biggl( \bigcap_{l+2 \leq i \leq n} E_i(z) \cap
\mathop{\bigcap
_{1 \leq j \leq n}}_{j \neq l-1,l,l+1} E_j(w) \biggr)\\
&&\qquad= \mathbf{P} \biggl( \bigcap_{l+2 \leq i \leq n} E_i(z) \biggr)
\mathbf{P} \biggl(
\mathop{\bigcap_{1 \leq j \leq n}}_{j \neq l-1,l,l+1} E_j(w) \biggr)\\
&&\qquad\leq \frac{1}{C^l \gamma_l^{a} s_l^{a}} \mathbf{P} \biggl( \bigcap
_{1 \leq i
\leq n}E_{i}(z) \biggr) \mathbf{P} \biggl( \bigcap_{1 \leq j \leq n}
E_{j}(w) \biggr).
\end{eqnarray*}
The last step comes by multiplying the second to last expression by
\[
\frac{1}{C^l \gamma_l^{a} s_l^{a}}\mathbf{P} \biggl(\bigcap_{1 \leq
i \leq l+1}
E_{i}(z) \biggr)\mathbf{P} \biggl( \bigcap_{l-1 \leq j \leq l+1}E_{j}(w)
\biggr)
\geq1
\]
and then using that each of the collections of events $E_i(z)$, $1 \leq
i \leq n$ and $E_j(w)$, $1 \leq j \leq n$ are independent (though, of
course, not from each other).
The lemma now follows as, by Proposition \ref{prop::markov},
$F_{n+1}(z)$ is independent of $E_m(z)$ for $1 \leq m \leq n$ and
$\mathbf{P}
(F_{n+1}(z)) \geq c > 0$ and the same is also true for $w$ with $c$
uniform in $n,z,w$.
\end{pf}

For $\alpha\geq0$, let $\nu_\alpha$ denote the Hausdorff-$\alpha$ measure.
\begin{lemma}
\label{lem::lower_bound}
We have\vspace*{1pt} $\mathbf{P}[\nu_{2-a}(T^C(a;U)) = \infty] = 1$ for all $0 <
a \leq2$
and $\mathbf{P}[\nu_2(T^C(0;U)) = \nu_2(U)] = 1$. In particular,
$\mathbf{P}[\dim
_H(T^C(a;U)) \geq2-a] = 1$ and $\mathbf{P}[ |T^C(2;U)| = \infty] = 1$.
\end{lemma}
\begin{pf}
Assume $ 0 < a \leq2$. Let $M_n = |H \cap C_n|$ and, for $z_{nj} \in H
\cap C_n$, let $p_{nj} = \mathbf{P}(z_{nj} \in C_n(a))$. For each $n
\in\mathbf{N}$,
define a random measure $\tau_n$ on $H$ by
\[
\tau_n(A) = \int_A \sum_{i=1}^{M_n} p_{ni}^{-1} 1_{ C_n(a)}(z_{ni})
1_{S(z_{ni},s_n)}(z) \,dz \qquad\mbox{for } A \subseteq H.
\]
Observe $\mathbf{E}\tau_n(H) = 1$ and
\begin{eqnarray*}
\mathbf{E}(\tau_n(H))^2
&=& s_n^4 \sum_{i,j=1}^{M_n} p_{ni}^{-1} p_{nj}^{-1} \mathbf{P}\bigl(z_{ni}, z_{nj}
\in C_n(a)\bigr)\\
&\leq& s_n^4 |M_n| \sum_{l \geq1} \biggl(\frac{ s_l^2}{s_n^2} \biggr)
O(C^l \gamma_l^{-a} s_l^{-a})\\
&=& \sum_{l \geq1} O(C^l \gamma_l^{-a} s_l^{2-a})
< \infty.
\end{eqnarray*}
Let
\[
I_\alpha(\tau_n) = \int_{[0,1]^2} \int_{[0,1]^2} \frac{ d\tau_n(z_1)
\,d\tau_n(z_2)}{|z_1 - z_2|^\alpha}
\]
be the $\alpha$-energy measure of $\tau_n$. By a similar computation,
we have
\begin{eqnarray*}
\mathbf{E}I_{\alpha}(\tau_n)
&=& \sum_{i,j=1}^{M_n} p_{ni}^{-1} p_{nj}^{-1} \mathbf{P}\bigl(z_{ni},
z_{nj} \in
C_{n}(a)\bigr) \int_{S(z_{ni},s_n)} \int_{S(z_{nj}, s_n)} \frac{ dz_1
\,dz_2}{|z_1-z_2|^\alpha}\\
&\leq&\sum_{l \geq1} O(C^l \gamma_l^{-a} s_l^{2-a} s_{l+1}^{-\alpha}),
\end{eqnarray*}
hence $\mathbf{E}I_{2-a}(\tau_n) < \infty$ uniformly in $n$. This
implies that
there exists $d,b >0$ such that with
\[
G_n = \{ b \leq\tau_n(H) \leq b^{-1}, I_{2-a}(\tau_n) \leq d\} \quad\mbox{and}
\quad G = \limsup_{n} G_n
\]
we have
\[
\mathbf{P}(G) > 0.
\]
As $I_{2-a}$ is lower semi-continuous the set $\mathcal{M}_{2-a}(b,d)$ of
measures $\tau$ on $H$ such that $b \leq\tau(H) \leq b^{-1}$ and
$I_{2-a}(\tau) \leq d$ is compact with respect to weak convergence. For
each $\omega\in G$ there exists a sequence $(n_k)$ such that $\tau
_{n_k,\omega} \in\mathcal{M}_{2-a}(b,d)$ and hence has a weak limit
$\tau
\in
\mathcal{M}_{2-a}(b,d)$ which is a finite measure supported on
$P(a)(\omega)$
with positive mass and finite $(2-a)$-energy. Therefore,
\[
\mathbf{P}\bigl(\nu_{2-a}(P(a)) > 0\bigr) > 0.
\]
A simple application of the Hewitt--Savage zero-one law implies that
\[
\mathbf{P}[ \nu_{2-a}(T^C(a;U)) > 0] = 1,
\]
and, in particular, $\mathbf{P}[ \dim_H(T^C(a;U)) \geq2-a] = 1$ (see
\cite{DPRZ01}, Lemma 3.2 for a similar argument).

We will now show that in fact $\mathbf{P}[\nu_{2-a}(T^C(a;U)) =
\infty] = 1$.
Consider the covering $S(z_{n_i},s_n)$ of $H$ by $(n!)^2$ disjoint
squares. The Markov property implies that we can write
$F|S(z_{n_i},s_n) = F_{n_i} + H_{n_i}$, where the $F_{n_i}$ are
independent zero-boundary GFFs and $H_{n_i}$ is the harmonic extension
of $F|\partial S(z_{n_i},s_n)$ to $S(z_{n_i},s_n)$. It is not hard to
see that $H_{n_i}$ is negligible in the definition of a thick point,
hence the set $T_{n_i}(a)$ of $a$-thick points of $F_{n_i}$ in
$S(z_{n_i},s_n)$ is the same as $T(a;U) \cap S(z_{n_i},s_n)$.
Therefore, the random variables $\nu_{2-a}(T_{n_i}(a))$ are i.i.d. and
$\nu_{2-a}(H) = \sum_i \nu_{2-a}(T_{n_i}(a))$. By the basic scaling
properties of $\nu_{2-a}$, we have that
\[
\nu_{2-a}(T_{n_i}(a)) \stackrel{d}{=} (s_n^{2-a}) \nu_{2-a}(H).
\]
The statement of the lemma in the case that $a > 0$ is now immediate.

It is left to consider the case that $a = 0$. It is immediate that
$\mathbf{P}
[z \in T^C(0;U)] = 1$ for any nonrandom $z \in U$. Hence, by Fubini's
theorem, we have that
\begin{eqnarray*}
\mathbf{E}\nu_2(T^C(0;U))
&=& \mathbf{E}\int_U 1_{T^C(0;U)}(z) \,d \nu_2(z)\\
&=& \int_U \mathbf{P}[z \in T^C(0;U)] \,d\nu_2(z) = \nu_2(U).
\end{eqnarray*}
Combining this with the trivial bound $0 \leq\nu_2(T^C(0;U)) \leq\nu
_2(U)$ implies\break $\nu_2(T^C(0;U)) = \nu_2(U)$ almost surely.
\end{pf}

\section{Conformal invariance}
\label{sec::conf_invariance}

The purpose of this section is to establish Theorem~\ref
{thm::conf_invariance} and Corollary \ref{cor::conf_invariance}. The
idea of the proof is to show that $\mu(A)$ is sufficiently well
approximated by $\sum_{n=1}^N \alpha_n f_n(z) |A|$ where $(f_n)$ is an
ONB of $H_0^1(U)$ and $A$ is a disk, square, or the conformal image of
such. The proof is divided into two subsections. In the first
subsection, we will compute the asymptotic variance of the GFF $F =
F_{[0,1]^2}$ on $[0,1]^2$ integrated over small disks, squares and the
conformal images of small disks and squares. We will then combine these
estimates with a covering argument and the Borel--Cantelli lemma to
bound $\mu(A) - \sum_{n=1}^N \alpha_n f_n(z) |A|$. The\vspace*{2pt} reason that we
restrict our attention to this case is that the $H_0^1([0,1]^2)$
orthonormal basis given by the eigenvectors of the Laplacian is
particularly convenient with which to work. In the second subsection,
we will combine these with some Gaussian estimates to prove the theorem.

\subsection{Preliminary estimates}
\label{subsec::estimates}

Let\vspace*{1pt} $F = F_{[0,1]^2}$ and let $\mu= \mu_{[0,1]^2}$ be given by $\mu(A)
= \int_{A} F(x)\,dx$. Throughout this section, we consider a fixed simply
connected domain $U$ and let $\varphi\dvtx U \to[0,1]^2$ be a
conformal transformation with inverse $\psi\dvtx[0,1]^2 \to U$. Fix
compact sets $K \subseteq U$ and $L \subseteq(0,1)^2$.
Let
\[
G_{ij}(z,r) = \int_{D(z,r)}\sin(\pi i u) \sin(\pi j v) \,du \,dv\qquad
\mbox{for } z \in[0,1]^2,
\]
and denote by $S(z,r)$ the square in $[0,1]^2$ centered at $z$ with
side length $r$.
\begin{lemma}
\label{lem::asymp_variance}
Uniformly in $z \in L$ and as $r \to0$,
%
\begin{eqnarray}\label{eqn::var_disk}
\mathbf{E}(\mu(D(z,r)))^2 &\sim&\frac{\pi}{2} r^4 \log\frac
{1}{r} ,\\
\label{eqn::var_square}
\mathbf{E}(\mu(S(z,r)))^2 &\sim&\frac{1}{2\pi} r^4 \log\frac{1}{r}.
\end{eqnarray}
\end{lemma}

We remark that it is possible to give a short proof of (\ref
{eqn::var_disk}) using (\ref{eqn::disk_integral}) and a little bit of
stochastic calculus. We give the following proof, however, because it
easily generalizes to the case of (\ref{eqn::var_square}) and the
intermediate estimates will be important for us later on.
\begin{pf*}{Proof of Lemma \protect\ref{lem::asymp_variance}}
We will only prove (\ref{eqn::var_disk}) as (\ref{eqn::var_square})
follows from the same argument. Using the representation (\ref
{eqn::gff_square_rep}), observe
\[
\mathbf{E}(\mu(D(z,r)))^2 = \frac{4}{\pi^2} \sum_{i,j \geq1}
\frac
{1}{i^2 +
j^2} G_{ij}^2(z,r).
\]
Let $g(r) = ( r \log\log\frac{1}{r} )^{-1}$,
\[
\Sigma_1 = \sum_{i,j \leq g(r)} \frac{1}{i^2 + j^2} G_{ij}^2(z,r)
\quad\mbox{and}\quad
\Sigma_2 = \sum_{i \vee j > g(r)} \frac{1}{i^2 + j^2} G_{ij}^2(z,r).
\]
With $z = (x,y)$ the symmetry of $D(0,r)$ implies
%
\begin{eqnarray} \label{eqn::front_estimate}\quad
G_{ij}(z,r)
&=& \int_{D(0,r)} \sin\bigl(\pi i (u+x)\bigr) \sin\bigl(\pi j (v+y)\bigr) \,du \,dv \nonumber
\\
&=& \int_{D(0,r)} \sin(\pi i x) \sin(\pi j y) \cos(\pi i u) \cos
(\pi j
v) \,du \,dv \nonumber\\[-8pt]\\[-8pt]
&=& \sin(\pi i x) \sin(\pi j y) \biggl( \pi r^2 + \int_{D(0,r)}[ \cos
(\pi
i u) \cos(\pi j v) - 1 ]\,du \,dv \biggr) \nonumber\\
&=& \sin(\pi i x) \sin(\pi j y) \bigl( \pi r^2 + O\bigl( r^4 (i \vee
j)^2\bigr)\bigr),\nonumber
\end{eqnarray}
so that for $i,j \leq g(r)$,
\[
G_{ij}^2(z,r) = \sin^2 (\pi i x) \sin^2 (\pi j y) \biggl( \pi^2 r^4 +
O\biggl(r^4 \biggl(\log\log\frac{1}{r}\biggr)^{-2}\biggr) \biggr) .
\]
Thus, by Lemma \ref{lem::sine_square_sum},
\[
\Sigma_1 \sim\frac{\pi^3}{8} r^4 \log\frac{1}{r}.
\]

We have
\[
|G_{ij}(z,r)|
= \biggl|\int_{x-r}^{x+r} \int_{y-\sqrt{r^2-x^2}}^{y+\sqrt{r^2- x^2}}
\sin(\pi i u) \sin(\pi j v) \,du \,dv \biggr|
\leq\int_{x-r}^{x+r} \frac{2}{\pi i}\,dv
= \frac{4 r}{\pi i }.
\]
Similarly, $|G_{ij}(z,r)| \leq\frac{4r}{\pi j}$ so that $G_{ij}(z,r) =
O ( \frac{r}{i \vee j} )$.
As
\[
\sum_{i \geq1} \sum_{j \geq g(r)} \frac{r^2}{(i^2 + j^2)(i \vee j)^2}
= O \biggl( \int_1^\infty\int_{g(r)}^\infty\frac{r^2}{(u^2 + v^2)u^2}
\,du \,dv \biggr)
\]
and
%
\begin{eqnarray} \label{eqn::tail_estimate}
\int_1^\infty\int_{g(r)}^\infty\frac{r^2}{(u^2 + v^2)u^2} \,du \,dv
&=& r^2 \int_{g(r)}^\infty\frac{1}{u^4} \int_1^\infty\frac
{1}{1+v^2/u^2} \,dv \,du\nonumber\\[-8pt]\\[-8pt]
&=& O\biggl(r^{4} \biggl(\log\log\frac{1}{r}\biggr)^{2}\biggr)\nonumber
\end{eqnarray}
it follows that $\Sigma_2$ is negligible compared to $\Sigma_1$ as $r
\to0$.
Therefore,
\[
\mathbf{E}(\mu(D(z,r)))^2 \sim\frac{4}{\pi^2} \Sigma_1 \sim\frac
{\pi
}{2} r^4
\log\frac{1}{r}.
\]
\upqed\end{pf*}

The purpose of the next lemma is to show that the same estimates hold
for \textit{conformal images} of disks and squares, the proof by simple
Fourier analysis. We will need to introduce some more notation. For
$\xi\in U$, let $\rho(r) = \rho(\xi,r) = |\varphi'(\xi)|r$,
$E(\xi,r)
= \varphi(D(\xi,r))$, $T(\xi,r) = \varphi(S(\xi,r))$, and
\[
H_{ij}(\xi,r) = \int_{E(\xi,r)} \sin(\pi i u) \sin(\pi j v) \,du \,dv.
\]
In the case of the former, we will always write $\rho(r)$ since $\xi$
will be clear from the context.
Obviously, the collection of functions
\[
(x,y) \mapsto2 \sin(\pi i x) \sin(\pi j y)
\]
is orthonormal in $L^2([0,1]^2)$ so that with $z = \varphi(\xi)$ Lemma
\ref{lem::conf_area_error} gives the bound
%
\begin{equation}
\label{eqn::l2_g_h_diff}
\sum_{i,j \geq1} \bigl(G_{ij}(z,\rho(r)) - H_{ij}(\xi,r)\bigr)^2 \leq C \int
_{[0,1]^2} \bigl|1_{D(z,\rho(r))} - 1_{E(\xi,r)}\bigr|^2
= O(r^3),\hspace*{-32pt}
\end{equation}
which holds uniformly in $\xi\in K$. As another consequence of Lemma
\ref{lem::conf_area_error}, we have
%
\begin{equation}
\label{eqn::g_h_diff}
|G_{ij}(z,\rho(r)) - H_{ij}(\xi,r)| = O(r^3) \qquad\mbox{as } r \to0,
\end{equation}
uniformly in $i,j$ and $\xi\in K$.
\begin{lemma}
\label{lem::asymp_variance_conf}
Uniformly in $\xi\in K$ we have
%
\begin{eqnarray} \label{eqn::im_var_disk}
\mathbf{E}(\mu(E(\xi,r)))^2 &\sim&\frac{\pi}{2} \rho^4(r) \log
\frac{1}{\rho(r)},\\
\label{eqn::im_var_square}
\mathbf{E}(\mu(T(\xi,r)))^2 &\sim&\frac{1}{2\pi} \rho^4(r) \log
\frac{1}{\rho(r)}.
\end{eqnarray}
\end{lemma}
\begin{pf}
As in the proof of Lemma \ref{lem::asymp_variance}, we will only show
(\ref{eqn::im_var_disk}) since the justification of (\ref
{eqn::im_var_square}) is exactly the same. Fix $\xi\in K$. For $z
=\varphi(\xi)$, let
\[
\Gamma_1 = \sum_{i,j \leq g(r)} \frac{G_{ij}^2(z,\rho(r)) -
H_{ij}^2(\xi,r)}{i^2 + j^2}
\]
and
\[
\Gamma_2 = \sum_{i \vee j>
g(r)} \frac{G_{ij}^2(z,\rho(r)) - H_{ij}^2(\xi,r)}{i^2 + j^2}.
\]
Using
\[
\sum_{i,j \leq n} \frac{1}{i^2 + j^2} =O(\log n),
\]
we see that (\ref{eqn::g_h_diff}) implies
%
\begin{eqnarray} \label{eqn::conf_front_estimate}\hspace*{28pt}
\Gamma_1
&=& \sum_{i,j \leq g(r)} \frac{1}{i^2 + j^2} \bigl(G_{ij}(z,\rho(r)) +
H_{ij}(\xi,r)\bigr)\bigl(G_{ij}(z,\rho(r)) - H_{ij}(\xi,r)\bigr)
\nonumber\\[-8pt]\\[-8pt]
&=& O(\log g(r) ) O(r^2) O(r^3)
= O\biggl(r^5 \log\frac{1}{r}\biggr).\nonumber
\end{eqnarray}
An application of (\ref{eqn::l2_g_h_diff}) and the Cauchy--Schwarz
inequality gives
\begin{eqnarray*}
&& \biggl|\sum_{i \vee j > g(r)} \bigl(G_{ij}^2(z,\rho(r)) - H_{ij}^2(\xi
,r)\bigr) \biggr|\\
&&\qquad\leq \biggl(\sum_{i,j \geq1} \bigl(G_{ij}(z,\rho(r)) + H_{ij}(\xi,r)\bigr)^2
\biggr)^{1/2}\\
&&\qquad\quad{}\times \biggl( \sum_{i,j \geq1} \bigl(G_{ij}(z,\rho(r)) -
H_{ij}(\xi
,r)\bigr)^2 \biggr)^{1/2} \nonumber\\
&&\qquad= [O(r^2) O(r^3)]^{1/2}
= O(r^{5/2}).
\end{eqnarray*}
Hence,
%
\begin{eqnarray} \label{eqn::conf_tail_estimate}
|\Gamma_2|
&\leq&\biggl(\sup_{i \vee j > g(r)} \frac{1}{i^2 + j^2} \biggr)
\biggl|\sum_{i,j > g(r)} \bigl(G_{ij}^2(z,\rho(r)) - H_{ij}^2(\xi,r)\bigr) \biggr|
\nonumber\\[-8pt]\\[-8pt]
& = & O\biggl(r^4 \biggl(\log\log\frac{1}{r}\biggr)^2\biggr).\nonumber
\end{eqnarray}
Therefore, uniformly in $\xi\in K$ and with $z = \varphi(\xi)$,
\[
\frac{1}{\rho^4(r) \log 1/{\rho(r)}} |\mathbf{E}
(\mu
(D(z,\rho
(r))) )^2 - \mathbf{E} (\mu(E(\xi,r)) )^2 | \to0
\qquad\mbox{as } r \to0.\hspace*{20pt}
\]
\upqed\end{pf}

Let $(\alpha_{ij})$ be the coefficients of $F$ as in (\ref
{eqn::gff_square_rep}) expressed in terms of the $H_0^1([0,1]^2)$
eigenbasis of $\Delta$. Let $r_n = e^{-n}$. For $r_{n+1} < r \leq r_n$,
set $\zeta(r) = r_{n}$, $\widehat{g}(r) = g(\zeta(r))$, and define
\[
\nu(A) = \mu(A) - \sum_{i,j \leq\widehat{g}(r)} \frac{2 \alpha
_{ij}}{\pi
\sqrt{i^2 + j^2}} |A| \sin(\pi i x) \sin(\pi j y),
\]
where $A$ is either a disk or a square centered at $z = (x,y) \in L$ of
radius $r$. If $A$ is the image of a disk or square centered at $\xi=
\psi(z) \in K$ of radius $r$ under $\varphi$, then we set
\[
\nu(A) = \mu(A) - \sum_{i,j \leq\widehat{g}(\rho(r))} \frac{2
\alpha
_{ij}}{\pi\sqrt{i^2 + j^2}} |A| \sin(\pi i x) \sin(\pi j y).
\]
The estimates (\ref{eqn::front_estimate}), (\ref{eqn::tail_estimate})
and (\ref{eqn::conf_front_estimate}), (\ref{eqn::conf_tail_estimate}) imply
\[
\mathbf{E}(\nu(A))^2 = O \biggl( r^4 \biggl(\log\log\frac{1}{r}\biggr)^2
\biggr).
\]
\begin{lemma}
\label{lem::asymp_gff_scaling}
Let $A$ be either a disk or square in $[0,1]^2$ centered in $L$ or the
image of such in $U$ under $\varphi$ centered in $K$. Then there exists
$\alpha= \alpha(\omega) > 0$ such that almost surely $\operatorname
{diam}(A) \leq
\alpha$ implies uniformly
%
\begin{eqnarray}
\label{eqn::mu_bound}
|\mu(A)| &=& O \biggl( |A| \log\frac{1}{|A|} \biggr),\\
\label{eqn::nu_bound}
|\nu(A)| &=& o \biggl( |A| \log\frac{1}{|A|} \biggr).
\end{eqnarray}
\end{lemma}

We will not make use of (\ref{eqn::mu_bound}) but record the result
anyway because its proof is the same as that of (\ref{eqn::nu_bound}).
\begin{pf*}{Proof of Lemma \protect\ref{lem::asymp_gff_scaling}}
We are going to give the complete proof in the case that $A$ is a disk
or square in $[0,1]^2$ centered in $L$ and then indicate the necessary
modifications to show that the result also holds for conformal images.
Lemma \ref{lem::asymp_variance} implies
\[
\mathbf{E}(\mu(S(z,2^{-n})))^2 \sim\frac{1}{2\pi} (2^{-n})^4 (\log2^n)
\]
so that for some $c_1 > 0$ and $n$ large enough
\[
\frac{ (2^{-n})^2 (\log2^n)}{\sqrt{\mathbf{E}(\mu
(S(z,2^{-n})))^2}} \geq c_1
\sqrt{n}.
\]
Therefore, by Lemma \ref{lem::norm_asymp_decay} with $c_2 = \sqrt{6}
c_1^{-1}$, we have
\[
\mathbf{P}\bigl(|\mu(S(z,2^{-n}))| > c_2 (2^{-n})^2 (\log2^n)\bigr)
= O( 2^{-3n}).
\]
Fix $\varepsilon> 0$ so that $L^\varepsilon$, the $\varepsilon
$-neighborhood of
$L$, satisfies $\overline{L^\varepsilon} \subseteq(0,1)^2$.
Letting $\mathcal{S}_n$ be the set of dyadic squares in $[0,1]^2$
contained in
$L^\varepsilon$ of side length $2^{-n}$ we see
\[
\sum_{n \geq1} \sum_{S \in\mathcal{S}_n} \mathbf{P}\bigl(|\mu(S)| >
c_2(2^{-n})^2
(\log
2^n)\bigr) < \infty.
\]
By the Borel--Cantelli lemma, there exists $n_0 = n_0(\omega)$ such
that for $n \geq n_0$ almost surely
%
\begin{equation}
\label{eqn::dyadic_square_bound}
|\mu(S)| \leq c_2 (2^{-n})^2 (\log2^n) \qquad\mbox{for all } S \in
\mathcal{S}_n.
\end{equation}

Suppose $R = [a,b] \times[c,d] \subseteq L^\varepsilon$ is a rectangle
with length $l = d-c$ and width $w= b-a = 2^{-n}$ with $n > n_0$ and $a
= i/2^n,b = (i+1)/2^n$ dyadic rationals. Assume further that $l \geq
w$. Fit as many dyadic squares of side length $2^{-n}$ into $R$ as
possible. Visibly, the number of such squares is bounded by $l/2^{-n}$.
The set that arises by removing these squares from $R$ consists of two
ends each of which contains at most
\[
\frac{2^{-n}}{2^{-n-1}} \cdot\frac{2^{-n}}{2^{-n-1}} = 2^{2}
\]
dyadic squares of side length $2^{-n-1}$. After removing these, each
end now contains at most
\[
\frac{2^{-n}}{2^{-n-2}} \cdot\frac{2^{-n-1}}{2^{-n-2}} = 2^{3}
\]
dyadic squares of side length $2^{-n-2}$. Iterating this procedure,
each end contains at most
\[
\frac{2^{-n}}{2^{-n-k}} \cdot\frac{2^{-n-(k-1)}}{2^{-n-k}} = 2^{k+1}
\]
squares of side length $2^{-n-k}$ at the $k$th step. Thus,
\begin{eqnarray*}
|\mu(R)|
&\leq& c_2 \biggl(\frac{l}{2^{-n}} (2^{-n})^2 (\log2^n) + 2\sum_{k
\geq
1} 2^{k+1} (2^{-n-k})^2 (\log2^{n+k}) \biggr)\\
&\leq& c_3 l w \log\frac{1}{w}.
\end{eqnarray*}
Now suppose that $2^{-n-1} < w \leq2^{-n}$ is not necessarily dyadic
and $l \geq2^{-n-1}$. Then a maximal decomposition of $R$ into
rectangles of length $l$ and with left- and right-hand sides located at
rationals of the form $i/2^{n+k}, (i+1)/2^{n+k}$, always taking the
largest possible such rectangle, contains at most two of width
$2^{-n-k}$ for each $k \in\mathbf{N}$ so that
%
\begin{eqnarray} \label{eqn::rectangle_bound_length}
|\mu(R)|
&\leq& c_4 l \biggl( \sum_{m \geq n+1} 2^{-m} \log2^{m} \biggr)
\leq c_5 l 2^{-n} n
\leq c_6 l w \log\frac{1}{l w}\nonumber\\[-8pt]\\[-8pt]
&=& c_6 |R| \log\frac{1}{|R|} \nonumber.
\end{eqnarray}
Note that this argument also works with the roles of $l$ and $w$ reversed.

Let $A$ be a disk contained in $L^\varepsilon$ with radius $r <
2^{-n_0-1}$. Slice $A$ vertically starting from the center to the right
and left into equal pieces of width $r^2$ and then slice it once
horizontally through the center. Let $A_1$ be the set consisting of the
union of the largest rectangles that fit into each slice. Then, since
there are at most $4r/r^2 = 4/r$ rectangles, each of area at most
$r^3$, (\ref{eqn::rectangle_bound_length}) gives us
\[
|\mu(A_1)| \leq12 c_6 r^2 \log\frac{1}{r}.
\]
Slice the regions above and below each of the rectangles in $A_1$,
including the degenerate rectangles on the left- and right-hand sides,
into equal pieces of width $r^3$. Denote the union of all the largest
rectangles contained in these slices by $A_2$ and note that the length
of each rectangle is at most $\sqrt{2r^3}$. The reason for this is that
the maximal length of such a rectangle with horizontal coordinates
contained in the interval $[a,b]$, say with $a \geq0$, is given by
$f(a) - f(b)$ where $f(x) = \sqrt{r^2 - x^2}$.
Obviously,
\[
f(a) - f(b) \leq f\bigl(r - (b-a)\bigr) - f(r) = f\bigl(r-(b-a)\bigr).
\]
In our case $b-a = r^2$ so that we have the bound $f(r - r^2) \leq
\sqrt
{ 2r^3}$.
If we iterate this procedure so that at the $n$th step we slice out
rectangles of width $r^{n+1}$ then at most $4 \cdot\frac{r}{r^{n+1}} =
4 r^{-n}$
rectangles each with length at most $f(r - r^{n}) \leq\sqrt{2} r^{(n+1)/2}$
and hence with area at most $\sqrt{2} r^{(n+1)/2} \cdot r^{n+1} =
\sqrt
{2} r^{3(n+1)/2}$.
If $A_n$ denotes the region from the $n$th step for $n \geq2$, then
\begin{eqnarray*}
|\mu(A_n)|
&\leq& 4\sqrt{2} c_6 r^{- n} \cdot r^{3(n+1)/2} \log\frac{1}{\sqrt
{2}r^{3(n+1)/2}}\\
&\leq& 12(n+1) c_6 r^{(n+3)/2} \log\frac{1}{r}.
\end{eqnarray*}
Therefore,
\[
|\mu(A)|
= \biggl|\mu\biggl(\bigcup_n A_n\biggr)\biggr|
\leq c_7 r^{3/2} \log\frac{1}{r} \sum_{n \geq1} n r^{n/2}
\leq c_{8} r^{2} \log\frac{1}{r}.
\]

This completes the proof of (\ref{eqn::mu_bound}) when $A$ is either a
disk or square centered in~$L$.

If $A$ is a square centered in $L$ of radius $r$, we know that
\[
\mathbf{E}(\nu(A))^2 = O \biggl( r^4\biggl(\log\log\frac{1}{r}\biggr)^2
\biggr) \qquad\mbox{as } r \to0,
\]
so that for some $d_1 > 0$ and $n$ large enough,
\[
\frac{ (2^{-n})^2 (\log2^n) (\log\log2^n)^{-1} }{\sqrt{ \mathbf
{E}(\nu
(S(z,2^{-n})))^2} } \geq d_1 \frac{n}{(\log n)^{2}} \geq d_1 \sqrt{n}.
\]
Hence for $d_2 > 0$ appropriately chosen and $n$ large enough,
\[
\mathbf{P}\bigl( |\nu(S(z,2^{-n}))| \geq d_2 (2^{-n})^2 (\log2^n) (\log
\log
2^n)^{-1}\bigr) = O(2^{-3n}).
\]
With $a(r) = (\log\log\frac{1}{r})^{-1}$
it follows from the Borel--Cantelli lemma that on dyadic squares small
enough and contained in $L^\varepsilon$ we have
\[
|\nu(S)| \leq d_2 |S| \biggl(\log\frac{1}{|S|} \biggr) a(|S|).
\]
If we do the covering argument as before, we can bound from above the
$\nu$-mass of the intermediate dyadic squares $S$ in our cover of $A$ by
\[
|\nu(S)| \leq d_2 |S| \biggl(\log\frac{1}{|S|} \biggr) a(|S|) \leq d_2
|S| \biggl(\log\frac{1}{|S|} \biggr) a(|A|).
\]
Thus, (\ref{eqn::nu_bound}) is now obvious.

To deduce the case when $A$ is a conformal image one runs the same
argument except instead of building coverings by dyadic squares in
$[0,1]^2$ one works with coverings by \textit{conformal images} of dyadic
squares in $U$. Indeed, we know by Lemma \ref{lem::asymp_variance_conf}
that the images of squares satisfy the same asymptotic variance bounds
as those in $L^\varepsilon$ up to a factor of $|\varphi'(\xi)|^2$. Hence,
one only needs to keep uniform control on $|\varphi'(\xi)|$ which is
easily accomplished by restricting to dyadic squares contained in a
neighborhood $K^\delta$ of $K$ such that $\overline{K^\delta}
\subseteq U$.
\end{pf*}

\subsection{Proof of conformal invariance}
\label{subsec::proof_of_conf}
Let $U,V \subseteq\mathbf{C}$ be bounded domains with smooth boundary and
$\varphi\dvtx U \to V$ a conformal transformation with inverse $\psi$.
\begin{pf*}{Proof of Theorem \protect\ref{thm::conf_invariance}}
Let $(S_n)$, $S_n = S(z_n,r_n)$, be a covering of $V$ by closed squares
such that $S(z_n,2r_n) \subseteq V$. Fix $K \subseteq U$ compact. With
$R_n = \psi(S_n)$, we can find indices $i_1,\ldots,i_k$ such that $K
\subseteq\bigcup_{1 \leq j \leq k} R_{i_j}$. Therefore, it suffices to show
%
\begin{equation}\qquad\quad
\label{eqn::conf_inv_restr}
\lim_{r \to0} \sup_{\xi\in K \cap R_{i_j}} \frac{1}{h(r)}
\bigl|\mu
_{U}(D(\xi,r)) - \mu_V(D(\varphi(\xi),|\varphi'(\xi)|r))|\psi
'(z)|^2 \bigr| = 0
\end{equation}
for each $j$ where $h(r) = \pi r^2 \log\frac{1}{r}$. If we write
$F_U|R_{i_j} = F_{i_j} + H_{i_j}$ with $F_{i_j}$ a zero-boundary GFF
and $H_{i_j}$ harmonic on $R_{i_j}$ then the term arising from
$H_{i_j}$ in (\ref{eqn::conf_inv_restr}) is negligible. As the same is
also true for $F_V|S_{i_j}$, we therefore may assume without loss of
generality that $U = \psi(S(z_{i_j},2r_{i_j}))$, which contains
$R_{i_j} \cap K$, and $V = S(z_{i_j},2r_{i_j})$. By a translation and
rescaling, we may further assume $V = [0,1]^2$.

For $\xi\in U$, let $E(\xi,r) = \varphi(D(\xi,r))$ be the image of the
disk $D(\xi,r) \subseteq U$ under~$\varphi$. With $\rho(r) =
|\varphi
'(\xi)|r$, Lemma \ref{lem::conf_area_error} implies $|E(\xi,r)
\Delta
D(\varphi(\xi),\rho(r))| = O(r^3)$ so that by Lemma \ref{lem::exp_norm}
we have
\begin{eqnarray*}
&&\mathbf{P} \biggl(\sup_{\xi\in K} \frac{1}{h(\rho(r))} \sum_{1
\leq i,j
\leq
g(r)} \frac{|\alpha_{ij}||E(\xi,r) \Delta D(\varphi(\xi),\rho(r))|}{
\sqrt{i^2 + j^2}} \geq t \biggr)\\
&&\qquad= \mathbf{P} \biggl( \sum_{1 \leq i,j \leq g(r)} \frac{|\alpha
_{ij}| O(r)}{
\sqrt
{i^2 + j^2}} \geq t \log\frac{1}{r} \biggr)\\
&&\qquad\leq r^t e^{O(r^2 \log g(r))} \prod_{1 \leq i,j \leq g(r)} \biggl(1+
\frac{O(r)}{\sqrt{i^2 + j^2}} \biggr).
\end{eqnarray*}
The inequality $\log(1+x) \leq x$ yields
\begin{eqnarray*}
\log\prod_{1 \leq i,j \leq g(r)} \biggl(1+ \frac{O(r)}{\sqrt{i^2 +
j^2}} \biggr) &\leq& O(r) \sum_{1 \leq i,j \leq g(r)} \frac{1}{\sqrt{i^2
+ j^2}} \\
&=& O\biggl(\biggl(\log\log\frac{1}{r}\biggr)^{-1}\biggr).
\end{eqnarray*}

Taking $r_n = e^{-n}$ and $t_n = n^{-1/2}$ it thus follows from the
Borel--Cantelli lemma that there exists $n_0 = n_0(\omega)$ such that
almost surely $n \geq n_0$ implies
\[
\sup_{\xi\in K} \frac{1}{h(\rho(r_n))} \sum_{1 \leq i,j \leq g(r_n)}
\frac{|\alpha_{ij}||E(\xi,r_n) \Delta D(\varphi(\xi),\rho
(r_n))|}{ \sqrt
{i^2 + j^2}} \leq\frac{1}{\sqrt{n}}.
\]
Combining this with Lemma \ref{lem::asymp_gff_scaling} yields for all
$n$ sufficiently large and $r_{n+1} < \rho(r) \leq r_n$ that
\begin{eqnarray*}
&& |\mu(E(\xi,r)) - \mu(D(z,\rho(r)))|\\[0.65pt]
&&\qquad \leq \sum_{1 \leq i,j \leq g(r_n)} \frac{2|\alpha_{ij}| |
|E(\xi
,r)| - |D(z,\rho(r))| |}{\pi\sqrt{i^2 + j^2}}\\[0.65pt]
&&\qquad\quad{} + |\nu(E(\xi
,r))| +
|\nu(D(z,\rho(r)))|\\[0.65pt]
&&\qquad\leq C \sum_{1 \leq i,j \leq g(r_n)} \frac{|\alpha_{ij}| |E(\xi,r_n)
\Delta D(z,\rho(r_n))| }{ \sqrt{i^2 + j^2}} + o(h(r))\\[0.65pt]
&&\qquad= o(h(r))
\end{eqnarray*}
uniformly in $\xi\in K$.
Therefore,
%
\begin{equation}
\label{eqn::uniform_conf_image_error}
\lim_{r \to0} \sup_{\xi\in K} \frac{1}{h(\rho(r))} |\mu
(E(\xi,r))
- \mu(D(\varphi(\xi),\rho(r))) | = 0.
\end{equation}
With $F_U = F \circ\varphi^{-1}$ the GFF on $U$ and $\mu_U(A) = \int_A
F_U(x) \,dx$, a change of variables gives
\[
\mu_U(D(\xi,r))
= \int_{[0,1]^2} F 1_{E(\xi,r)} | \psi'|^2
= [|\psi'(z)|^2 + O(r)] \mu(E(\xi,r))
\]
uniformly in $\xi\in K$. With $z = \varphi(\xi)$ we have
\begin{eqnarray*}
&& \bigl|\mu_U(D(\xi,r)) - \mu(D(z,\rho(r)))|\varphi'(\xi)|^{-2}\bigr|\\[0.65pt]
&&\qquad= \bigl| \mu(E(\xi,r))[ |\psi'(z)|^2 + O(r)] - \mu(D(z,\rho
(r)))|\varphi
'(\xi)|^{-2}\bigr|\\[0.65pt]
&&\qquad\leq \bigl| \mu(E(\xi,r)) |\psi'(z)|^2 - \mu(D(z,\rho(r)))|\varphi
'(\xi
)|^{-2}\bigr| + O(r)|\mu(E(\xi,r))|.
\end{eqnarray*}
The theorem now follows as by (\ref{eqn::uniform_conf_image_error}),
we have
\begin{eqnarray*}
\mu(E(\xi,r)) |\psi'(z)|^2
&=& \mu(D(z,\rho(r)))|\psi'(z)|^2 + o(h(r))\\[0.65pt]
&=& \mu(D(z,|\varphi'(\xi)|r))|\varphi'(\xi)|^{-2}
+ o(h(r)).
\end{eqnarray*}
\upqed\end{pf*}

\begin{appendix}
\section{Gaussian estimates}

\begin{lemma}
\label{lem::sine_square_sum}
If $L \subseteq(0,1)^2$ is compact, then
\[
\sum_{1 \leq i,j \leq n} \frac{\sin^2 (\pi i x) \sin^2 (\pi
jy)}{i^2 +
j^2} \sim\frac{\pi}{8} \log n \qquad\mbox{as } n \to\infty,
\]
uniformly in $(x,y) \in L$.
\end{lemma}
\begin{pf}
Observe there exists $c > 0$ such that
\[
\sum_{1 \leq i,j \leq n} \frac{\sin^2 (\pi i x) \sin^2 (\pi j
y)}{i^2 +
j^2} \geq c \log n.
\]
Hence, as far as the asymptotics of the summation are concerned, we may
ignore terms that are $o(\log n)$. Thus as $n \to\infty$, we have
\begin{eqnarray*}
&&\sum_{1 \leq i,j \leq n} \frac{\sin^2 (\pi i x) \sin^2 (\pi j y)}{i^2
+ j^2}
\\
&&\qquad\sim\int_1^n \int_1^n \frac{\sin^2 (\pi u x) \sin^2 (\pi v
y)}{u^2 +
v^2} \,du \,dv\\
&&\qquad\sim \sum_{1 \leq i,j \leq n} \int_{{i}/{x}}^{
({i+1})/{x}} \int
_{{j}/{y}}^{({j+1})/{y}} \frac{\sin^2 (\pi u x) \sin^2
(\pi v
y)}{u^2 + v^2} \,du \,dv\\
&&\qquad\sim \sum_{1 \leq i,j \leq n} \frac{1}{(i/x)^2 + (j/y)^2} \int
_{{i}/{x}}^{({i+1})/{x}} \int_{{j}/{y}}^{({j+1})/{y}} \sin^2
(\pi u x) \sin^2 (\pi v y) \,du \,dv\\
&&\qquad= \frac{1}{4} \sum_{1 \leq i,j \leq n} \frac{1}{(i/x)^2 + (j/y)^2}
\frac{1}{x} \cdot\frac{1}{y}\\
&&\qquad\sim \frac{1}{4} \sum_{1 \leq i,j \leq n} \int_{
{i}/{x}}^{
({i+1})/{x}} \int_{{j}/{y}}^{({j+1})/{y}} \frac{1}{u^2 + v^2}
\,du \,dv\\
&&\qquad\sim\frac{1}{4} \int_1^n \int_1^n \frac{1}{u^2 + v^2} \,du \,dv\\
&&\qquad\sim \frac{\pi}{8} \log n.
\end{eqnarray*}
\upqed\end{pf}
\begin{lemma}
\label{lem::exp_norm}
If $(X_n)$ is an i.i.d. sequence of standard normals and $(\beta_n)$ a
sequence of positive constants, then we have
\[
\mathbf{P} \biggl( \sum_{n} \beta_n |X_n| \geq t \biggr)
\leq e^{-t} \prod_n (1+\beta_n)e^{\beta_n^2/2}.
\]
\end{lemma}
\begin{pf}
Markov's inequality gives
\[
\mathbf{P} \biggl( \sum_{n} \beta_n |X_n| \geq t \biggr)
\leq e^{-t} \prod_n \mathbf{E}\exp( \beta_n |X_n| ).
\]
If $X \sim N(0,1)$ and $\beta> 0$, we have
\[
\mathbf{E}e^{\beta X} 1_{\{X \geq0\}}
= \frac{e^{\beta^2/2}}{\sqrt{2\pi}} \int_{-\beta}^\infty e^{-
x^2/2} \,dx
\leq\biggl(\frac{1}{2} + \frac{\beta}{\sqrt{2\pi}} \biggr)
e^{\beta^2/2}.
\]
Thus, $\mathbf{E}e^{\beta|X|} \leq(1+\beta) e^{\beta^2/2}$. Combining
everything gives the result.
\end{pf}
\begin{lemma}
\label{lem::bm_bound}
Let $B(t)$ be a standard Brownian motion, $\mu> 0$, and $T \geq1$
fixed. Then
\[
\mathbf{P}\Bigl( {\max_{0 \leq t \leq T}} |B(t) - \mu t| \leq\sqrt{T}\Bigr)
\geq C \exp
\biggl(\frac{1}{2}\bigl(\mu\sqrt{T}-\mu^2 T\bigr)\biggr),
\]
where $C > 0$ is a constant independent of $T$.
\end{lemma}
\begin{pf}
Let $E_T^\mu= \{ {\max_{0 \leq t \leq T}} |B(t) - \mu t| \leq\sqrt{T}
\}$.
By the Girsanov theorem,
\[
\mathbf{P}(E_T^\mu)
= \mathbf{E}\bigl[ e^{\mu B(T) - \mu^2T/2} 1_{E_T^0}\bigr]
\geq\exp\bigl(\tfrac{1}{2}\bigl(\mu\sqrt{T}-\mu^2 T \bigr)\bigr) \mathbf{P}\bigl[ E_T^0, B(T)
\geq
\sqrt{T}/2 \bigr].
\]
Taking
\[
C = \mathbf{P}\bigl[ E_T^0, B(T) \geq\sqrt{T}/2 \bigr] = \mathbf{P}\Bigl( {\max_{0
\leq t \leq1}}
|B(t)| \leq1, B(1) \geq1/2\Bigr) > 0
\]
proves the lemma.
\end{pf}
\begin{lemma}
\label{lem::norm_asymp_decay}
If $Z \sim N(0,1)$, then
\[
\mathbf{P}( |Z| > \lambda) \sim\sqrt{\frac{2}{\pi}} \lambda^{-1}
e^{-\lambda
^2/2} \qquad\mbox{as } \lambda\to\infty.
\]
\end{lemma}
\begin{pf}
See Lemma 1.1 of \cite{OP73}. 
\end{pf}

\section{Area distortion under conformal maps}

\begin{lemma}
\label{lem::conf_area_error}
Suppose that $U,V \subseteq\mathbf{C}$ are domains with $K \subseteq U$
compact. If $\varphi\dvtx U \to V$ is a conformal transformation, then
\[
| E(\xi,r) \Delta D(\varphi(\xi),\rho(r))| = O(r^3)
\]
uniformly in $\xi\in K$ where $E(\xi,r) = \varphi(D(\xi,r))$ and
$\rho
(r) = \rho(\xi,r) = |\varphi'(\xi)|r$.
\end{lemma}
\begin{pf}
For $|\xi-\eta| \leq r$, we have
\[
|\varphi(\xi) - \varphi(\eta)| = |\varphi'(\xi)(\xi-\eta) +
O(r)(\xi
-\eta)|
\]
so that
\[
\bigl(|\varphi'(\xi)| - O(r)\bigr)|\xi-\eta| \leq|\varphi(\xi) - \varphi
(\eta)|
\leq\bigl(|\varphi'(\xi)| + O(r)\bigr)|\xi-\eta|.
\]
This implies
\[
D\bigl(\varphi(\xi),\rho(r) - O(r^2)\bigr) \subseteq E(\xi,r) \subseteq
D\bigl(\varphi
(\xi),\rho(r)+O(r^2)\bigr),
\]
which gives
\begin{eqnarray*}
|E(\xi,r) \Delta D(\varphi(\xi),\rho(r))|
&\leq& \bigl|D\bigl(\xi,\rho(r)+O(r^2)\bigr)\bigr| - \bigl|D\bigl(\xi,\rho(r) - O(r^2)\bigr)\bigr|\\
&=& O(r^3) \qquad\mbox{as } r \to0
\end{eqnarray*}
uniformly in $z \in K$.\vspace*{-18pt}
\end{pf}

\section{Modified Kolmogorov--Centsov}

\begin{lemma}
\label{lem::kolm_centsov} Suppose that $U \subseteq\mathbf{R}^d$ is
a bounded
open set and that $X \dvtx U \times(0,1] \to\mathbf{R}$ is a time-varying
random field satisfying
\[
\mathbf{E}| X(z,r) - X(w,s)|^{\alpha} \leq C \biggl( \frac{|(z,r) - (w,s)|}{r
\wedge s} \biggr)^{d+1 + \beta}
\]
for some $\alpha, \beta> 0$. For each $\zeta> \alpha^{-1}$ and
$\gamma\in(0,\beta/\alpha)$, $X$ has a modification $\tilde{X}$ satisfying
\[
|\tilde{X}(z,r) - \tilde{X}(w,s)| \leq M \biggl( \log\frac{1}{r}
\biggr)^\zeta\frac{|(z,r) - (w,s)|^\gamma}{r^{\tilde{\gamma}}},
\]
where
\[
\tilde{\gamma} = \frac{d+\beta}{\alpha},
\]
$z,w \in U$ and $r,s \in(0,1]$ with $1/2 \leq r/s \leq2$.
\end{lemma}

The proof is almost exactly the same as the usual proof of the
Kolmogorov--Centsov theorem \cite{KS98,RY04}. We will include a
proof for the convenience of the reader which will follow very closely
that given in \cite{RY04}.
\begin{pf*}{Proof of Lemma \protect\ref{lem::kolm_centsov}}
We may assume without loss of generality that $U \subseteq[0,1]^d$ by
rescaling. For each $n,T \in\mathbf{N}$, let
\[
R_n^T = \{(\underline{i},j)/ 2^n \in U \times(2^{-T},2^{1-T}] \dvtx
\underline{i} \in\mathbf{Z}^d,j \in\mathbf{Z}\} \quad\mbox{and}\quad R^T
= \bigcup_n R_n^T.
\]
Let $\Delta_n^T$ be the set of pairs $a,b \in R_n^T$ such that $|a - b|
= 2^{-n}$. Trivially, $|\Delta_n^T| = O( 2^{(n+1)(d+1) - T})$. Let
\[
K_i = \sup_{T \geq1} \biggl( \frac{2^{-\tilde{\gamma} T}}{T^\zeta
} {\sup
_{(a,b) \in\Delta_i^T}}|X(a) - X(b)| \biggr).
\]
We have
\begin{eqnarray*}
\mathbf{E}K_i^\alpha
&\leq& \sum_{T \geq1} \frac{2^{- \alpha\tilde{\gamma}
T}}{T^{\alpha
\zeta}} \sum_{a,b \in\Delta_i^T} \mathbf{E}|X(a) - X(b)|^\alpha\\
&\leq& \sum_{T \geq1} \frac{C 2^{-(d+\beta) T}}{T^{\alpha\zeta}}
2^{(i+1)(d+1) - T} \cdot2^{(T-i)(d+1+\beta)}
= O(2^{-i\beta}).
\end{eqnarray*}
For $a,b \in U \times(0,1]$, we say that $a \leq b$ if the
corresponding component-wise inequalities hold. If $a \in R^T$, then
there exists an increasing sequence $(a_n)$ in $U \times
(2^{-T},2^{1-T}]$ such that $a_n \in R_n^T$ for every $n$, $a_n \leq
a$, and $a_n = a$ for all $n$ large enough. Let $b \in R^T$ and $(b_n)$
be such a sequence for $b$. Assume $|a-b| \leq2^{-m}$. Then
\begin{eqnarray*}
X(a) - X(b)
&=& \sum_{i=m}^\infty\bigl( X(a_{i+1}) - X(a_i)\bigr) + X(a_m) - X(b_m) \\
&&{} + \sum_{i=m}^\infty\bigl(X(b_i) - X(b_{i+1})\bigr),
\end{eqnarray*}
which implies
\[
\frac{2^{-\tilde{\gamma} T }}{T^\zeta} |X(a) - X(b)| \leq K_m + 2
\sum
_{i=m+1}^\infty K_i \leq2 \sum_{i=m}^\infty K_i.
\]
We have
\begin{eqnarray*}
A
&\equiv& \sup_{T,m} \biggl( \sup\biggl\{\frac{2^{- \tilde{\gamma}
T} |X(a)
- X(b)| }{T^{\zeta} |a - b|^\gamma} \dvtx a,b \in R^T, 2^{-(m+1)} \leq
|a-b| \leq2^{-m} \biggr\} \biggr)\\
&\leq& \sup_{m \in\mathbf{N}} \Biggl( 2^{\gamma(m+1) + 1} \sum
_{i=m}^\infty K_i
\Biggr)
\leq2 \sum_{i=0}^\infty2^{\gamma i} K_i.
\end{eqnarray*}
This implies $\mathbf{E}A^\alpha< \infty$ so that for some $M > 0$ when
$(z,r), (w,s) \in R^T$ we have
\[
|X(z,r) - X(w,s)|
\leq M T^\zeta\frac{|(z,r) - (w,s)|^\gamma}{2^{-\tilde{\gamma}T}}
\leq M \biggl( \log\frac{1}{r} \biggr)^{\zeta} \frac{|(z,r) -
(w,s)|^\gamma}{r^{\tilde{\gamma}}}.
\]
With $R = \bigcup_T R^T$,
\[
\tilde{X}(a) = \mathop{\lim_{b \to a}}_{ b \in R} X(b)
\]
is clearly is the desired modification.
\end{pf*}
\end{appendix}

\section*{Acknowledgments}
We thank Scott Sheffield for a helpful discussion that led to the
observation that the law of the Hausdorff-$(2-a)$ measure of the
$a$-thick points is infinitely divisible, which in turn led to the
proof that the Hausdorff-$(2-a)$ measure of the $a$-thick points is
infinite. We also thank Charles Newman for posing the question of
whether or not the set of $2$-thick points was empty.

%

%
\printaddresses

\end{document}